\documentstyle[12pt]{article}
\begin{document}
\author{S.V. Ludkovsky.}
\title{Stochastic processes on 
geometric loop groups and diffeomorphism groups 
of real and complex manifolds,
associated unitary representations.}
\date{27 January 2001}
\maketitle
\begin{abstract}
This article is devoted to the investigation of the 
Belopolskaya's and Dalecky's problem. It consists of construction 
of a stochastic process on an infinite dimensional Lie group
$G$ which does not satisfy locally the Campbell-Hausdorff formula
and construction of a dense subgroup $G'$ in $G$ such that a 
transition measure is quasi-invariant and differentiable 
relative to the left or right action of $G'$ on $G$.
Geometric loop groups and diffeomorphism groups 
of Sobolev classes of smoothness are investigated 
for finite dimensional and also infinite dimensional 
real manifolds. Such groups also are defined and studied for
complex manifolds finite and infinite dimensional.
Stochastic processes are considered on free loop spaces, 
geometric loop and diffeomorphism groups of real and complex manifolds.
They are used for investigations of Wiener
differentiable quasi-invariant measures on such groups 
relative to dense subgroups. 
Such measures are used for the investigation of associated unitary 
representations of these groups. 
\end{abstract}
\section{Introduction.}  
Earlier Gaussian quasi-invariant measures on loop groups
of Riemann and complex manifolds were investigated \cite{ludan,lulgcm,lulgrm}.
With the help of them irreducible strongly continuous unitary 
representations were constructed. Gaussian measures were studied 
on free loop spaces also. 
Traditionally geometric loop groups are considered
as families of mappings $f: S^1\to N$ from the unit circle
into a Riemann manifold $N$ preserving marked points $s_0\in S^1$ and 
$y_0\in N$ under the corresponding equivalence relation
caused by an action of a diffeomorphism group $Diff^{\infty }(S^1)$
of the circle on the free loop space \cite{gajer}. 
But in \cite{ludan,lulgcm,lulgrm} were defined and 
investigated generalized loop groups as families of mappings from one 
manifold $f: M\to N$ into another preserving marked points
$s_0\in M$ and $y_0\in N$ under the corresponding equivalence
relation in the free (pinned) loop space and with the help of Grothendieck 
construction of an Abelian group from a commutative monoid with the 
unit and the cancellation property with rather mild conditions on 
$M$ and $N$ for finite and infinite dimensional real and complex manifolds.
\par If consider the composition of two nontrivial pinned in the marked point
$s_0$ loops of class $C^n$,
where $n\ge 1$, then the resulting loop is continuous, 
but generally not of class  
$C^n$ as it can be lightly seen on examples of the unit circle
$S^1$ and the unit sphere $S^2$.
If for the $S^1$ case $f$ is a $C^n$-loop with $n\ge 1$, then $f$ and $f'$
are continuous functions by the polar coordinate $\theta \in [0,1]$
such that $f(0)=f(1)$ and $\lim_{\theta \to 0, \theta >0}f'(\theta )=
\lim{\theta \to 1, \theta <1}f'(\theta )=:f'(0)$. There are another
$C^n$-loops $g$ such that $g'(0)\ne f'(0)$, but $g(0)=f(0)$. Then a loop
$h(\theta ):=f(2\theta )$ for each $\theta \in [0,1/2)$ and
$h(\theta ):=g(2\theta -1)$ for each $\theta \in [1/2,1]$
is a continuous loop, but not a $C^1$-loop. 
This is generally only continuous and piecewise of class
$C^n$ and such submanifolds, restrictions on which
are of class $C^n$, can be described as submanifolds with corners.
Another reason is that $S^n\vee S^n$ is a continuous retraction
of $S^n$, but there is not any diffeomorphism between them.
Also from $S^1\times S^n$ there is a continuous mapping on
$S^{n+1}$, but it is not a diffeomorphism (see \S 2.1.4).
Naturally, for a definition
of the smooth composition in loop monoids and loop groups
manifolds with corners (with the corresponding atlases)
are used. This permits to define topological loop monoids and 
topological loop groups. Another two reasons of the consideration of 
manifolds with corners are given below.
\par The commutative monoid is not the free (pinned) loop space, since
it is obtained from the latter by factorization.
For the construction of loop groups here are used manifolds $M$
with some mild additional conditions.
When $M$ is finite dimensional over $\bf C$ we suppose that it is compact.
This condition is not very restrictive, since each locally compact
space has Alexandroff (one-point) compactification (see Theorem 3.5.11
in \cite{eng}).
When $M$ is infinite dimensional over $\bf C$ it is assumed,
that $M$ is embedded as a closed bounded subset into the 
corresponding Banach space $X_M$ over $\bf C$.
This is necessary that to define a group structure on a 
quotient space of a free loop space. 
\par The free loop space
is considered as consisting of continuous functions $f: M\to N$
which are (piecewise) 
holomorphic in the complex case or 
(piecewise) continuously differentiable in the real case
on $M\setminus M'$ and preserving marked points
$f(s_0)=y_0$, where $M'$
is a closed real submanifold depending on $f$ with a codimension
$codim_{\bf R}M'=1$, $s_0\in M$ and $y_0\in N$ are marked points. 
There are two reasons to consider
such class of mappings. The first is the need to define correctly
compositions of elements in the loop group (see beneath).
The second is the isoperimetric inequality for holomorphic
loops, which can cause the condition of a loop
to be constant on a sufficiently small neighbourhood of $s_0$
in $M$, if this loop is in some small neighbourhood of $w_0$,
where $w_0(M):=\{ y_0 \} $ is a constant loop 
(see Remark 3.2 in \cite{hum}).
\par In this article loop groups of different classes are
considered. Classes analogous to Gevrey classes
and also with the usage of Sobolev classes of $f: M\setminus M' \to N$  
are considered for the construction of dense loop 
subgroups and quasi-invariant measures. Henceforth, we consider
not only orientable manifolds $M$ and $N$, but also nonorientable
manifolds (apart from \cite{lulgcm}, where only orientable 
manifolds were considered), since 
for a non-orientable manifold there always exists
its orientable double covering manifold (see \S 6.5 
in \cite{abma}).
Loop commutative monoids with the cancellation property
are quotients of families of mappings $f$ from $M$ into a
manifold $N$ with $f(s_0)=y_0$ by the corresponding equivalence
relation. 
For the definition of the equivalence relation here are 
not used groups of holomorphic diffeomorphisms because of strong
restrictions on their structure caused by 
holomorphicity (see Theorems 1 and 2 in \cite{bomon}).
Groups are constructed from monoids with the help of A. Grothendieck
procedure.
These groups are commutative and non-locally compact.
They does not have non-trivial local one-parameter subgroups
$\{ g^b:$ $b\in (-a,a) \} $ with $a>0$ for an element $g$ corresponding
to a class of a mapping $f: M\to N$, $f(s_0)=y_0$, when 
$f$ is such that $\sup_{y\in N} [card(f^{-1}(y))]=k<\aleph _0$, since 
$g^{1/p}$ does not exist in the loop group for each prime integer
$p$ such that $p>k$ (see \S 2).
Therefore, in each neighbourhood
$W$ of the unit element $e$ there are elements which does not belong
to any local one-parameter subgroup.
\par These groups are Abelian, non-locally compact and 
for them the Campbell-Hausdorff formula is not valid (in 
an open local subgroup). Finite dimensional Lie groups 
satisfy locally the Campbell-Hausdorff formula.
This is guarantied, if impose on a locally compact topological
Hausdorff group $G$ two conditions: it is a $C^{\infty }$-manifold
and the following mapping $(f,g)\mapsto f\circ g^{-1}$
from $G\times G$ into $G$ is of class $C^{\infty }$.
But for infinite dimensional $G$ the Campbell-Hausdorff 
formula does not follow from these conditions.
Frequently topological Hausdorff groups satisfying these 
two conditions also are called Lie groups, though they 
can not have all properties of finite dimensional 
Lie groups, so that the Lie algebras for them do not play 
the same role as in the finite dimensional case and therefore
Lie algebras are not so helpful.
If $G$ is a Lie group and its tangent space $T_eG$ is a Banach space,
then it is called a Banach-Lie group, sometimes it is undermined,
that they satisfy the Campbell-Hausdorff formula locally for a Banach-Lie 
algebra $T_eG$. In some papers the Lie group terminology
undermines, that it is finite dimensional.
It is worthwhile to call Lie groups satisfying the Campbell-Hausdorff
formula locally (in an open local subgroup) by Lie groups in the 
narrow sense;
in the contrary case to call them by Lie groups in the broad sense.
\par Stochastic processes on Lie groups $G$ were considered in
\cite{ahmt,beldal,dalschn}. The book \cite{ahmt} is devoted also
to Lie algebras and to Lie groups satisfying the Campbell-Hausdorff 
formula, the theory of which differs drastically from the groups considered
in this paper. General theorems about quasi-invariance 
and differentiabilty of transition measures on the Lie group $G$
relative to a dense subgroup $G'$ were given in \cite{beldal,dalschn}, 
but they 
permit to find $G'$ only abstractly and when a local subgroup of $G$ 
satisfies the Campbell-Hausdorff formula. For Lie groups which 
do not satisfy the Campbell-Hausdorff formula locally
this question remained open, as was pointed out by Belopolskaya 
and Dalecky in Chapter 6. They have proposed in such cases to investigate
concrete Lie groups that to find pairs $G$ and $G'$.
As it is well-known in mathematics problems of an existence
of an object and a description of it are frequently called
perpendicular.
In each concrete case of  $G$ it its necessary to construct a stochastic 
process and $G'$. On the other hand, the groups considered in the present
article do not satisfy the Campbell-Hausdorff formula.
\par Below path spaces, loop spaces, pinned loop groupoids,
loop monoids, loop groups and diffeomorphism groups are considered
not only for finite dimensional, but also for infinite dimensional 
manifolds. The path spaces also are called path groups,
but they have the group structure neither in the usual algebraic
sence nor in the usual topological group sence. Path spaces are more or 
less known (see for example, \cite{aht,malb}). In this article 
they are mentioned mainly 
from the manifold point of view 
and in the generalized sense for the considered here classes of 
smoothness and manifold structures compatible with the 
manifold structures of loop groups.
\par In particular, loop and diffeomorphism groups are important for the 
development of the representation theory of non-locally compact groups. 
Their representation theory has many 
differences with the traditional representation theory of 
locally compact groups and finite dimensional Lie groups, 
because non-locally 
compact groups have not $C^*$-algebras associated with the Haar 
measures and they have not underlying Lie algebras and relations between 
representations of groups and underlying algebras (see also \cite{lubp}). 
\par In view of the A. Weil 
theorem if a topological Hausdorff group $G$ has a quasi-invariant measure 
relative to the entire $G$, then $G$ is locally compact.
Since loop groups $(L^MN)_{\xi }$ are not locally compact, they can not have 
quasi-invariant measures relative to the entire group, but only relative to 
proper subgroups $G'$ which can be chosen dense in  $(L^MN)_{\xi }$,
where an index $\xi $ indicates on a class of smoothness.
The same is true for diffeomorphism groups (besides holomorphic 
diffeomorphism groups of compact complex manifolds).
Diffeomorphism groups of compact complex manifolds
are finite dimensional Lie groups (see \cite{kobtg}
and references therein).
It is necessary to note that there are quite another groups
with the same name loop groups, but they are infinite dimensional
Banach-Lie groups of mappings $f: M\to H$ into a finite dimensional Lie group
$H$ with the pointwise group multiplication of mappings with values in $H$.
The loop groups considered here are generalized geometric loop groups.
\par The traditional geometric loop groups and free loop spaces are important both 
in mathematics and in modern physical theories. Moreover,
generalized geometric loop groups also can be used in the same fields
of sciences
and open new opportunities. In cohomology theory and physical applications
stochastic processes and Wiener measures on the free loop spaces are used
\cite{aiel,eljm,hsufa,jole,jolesa,lejs,lejmp,malb}.
In these papers were considered only particular cases of 
real free loop spaces and groups for finite dimensional manifolds,
no any applications to the representation theory were given.
\par On the other hand, representation theory of non-locally compact groups
is little developed apart from the case of locally compact groups.
For locally compact groups theory of induced representations
is well developed due to works of Frobenius, Mackey, etc.
(see \cite{barut,fell} and references therein).
But for non-locally compact groups it is very little known.
In particular geometric loop and diffeomoprphism groups have important 
applications in modern physical theories (see \cite{ish,mensk}
and references therein).
\par One of the main tools in the investigation of unitary represenations 
of nonlocally compact groups are quasi-invariant measures.
In previous works of the author \cite{lurim1,lurim2} 
Gaussian quasi-invariant 
measures were constructed on diffeomorphism groups
with some conditions on real manifolds. For example, compact manifolds 
without boundary were not considered, as well as infinite dimensional 
manifolds with boundary. In this article new Gevrey-Sobolev classes 
of smoothness for diffeomorphism groups of infinite dimensional 
real and complex manifolds 
are defined and investigated. This permits to define on them the Hilbert 
manifold structure. This in its turn simplifies the construction of 
stochastic processes and transition quasi-invariant measures on them.
Wiener transition quasi-invariant measures are constructed below for 
wider classes of manifolds. Pairs of topological groups $G$ and their dense
subgroups $G'$ are described precisely.
\par This work is devoted to the investigation of Wiener measures 
and stochastic processes on the generalized loop spaces,
loop monoids, geometric loop groups and diffeomorphism groups. 
For the loop groups are considered both measures
arising from the stochastic equations on them and aslo induced from 
the free loop space. Their quasi-invariance and differentiability 
relative to dense subgroups is investigated.
Transition measures arising from stochastic processes on
manifolds also are called Wiener measures.
Then measures are used for the study of associated unitary 
regular and induced representations of dense subgroups $G'$. 
\par Section 2 is devoted to the definitions of topological and 
manifold structures of loop groups and diffeomorphism groups
and their dense subgroups.
In section 3 Wiener processes and transition quasi-invariant
differentiable measures are studied
(see Theorem 3.3). 
Unitary representations of dense subgroups $G'$ 
founded in sections 2 and 3 are investigated in section 5.
\par Section 4 is devoted to loop monoids
as well as to loop groupoids, which are defined 
in \S 4.2. For the considered here classes of manifolds
the generalized path space is defined in \S 4.4.
All objects given in sections 2-4 were
not considered by others authors, besides very specific particular 
cases of the diffeomorphism group
$Diff^{\infty }(S^1)$ and loop groups for $M=S^1$ and path spaces
for $M=[0,1]$ outlined above.
Differentiable transition Wiener measures on them
are given in Theorems 4.1, 4.3 and 4.5.
Basic facts and notations of stochastic analysis on manifolds
are reminded in the Appendix, that may be useful,
for example, for specialists in group theory or differential geometry
do not working with stochastic analysis. 
\section{Loop and diffeomorphism groups of real and complex 
finite and infinite dimensional manifolds.}
To avoid misunderstandings we first give our definitions
of manifolds considered here and then of loop and diffeomorphism
groups. In \S 2.1.1 uniform atlases are defined.
They are necessary on Lie groups for the construction of 
stochastic processes on them. 
In \S \S 2.1.2-2.1.5 loop groups and in
\S \S 2.2-2.4 diffeomorphism groups are defined.
In \S \S 2.1.6-2.1.8, 2.5-2.9 necessary statements about
their structures as Lie groups and manifolds are given.
\par For loop groups and diffeomorphism groups manifolds 
are supposed to be satisfying the corresponding specific conditions.
They are related mainly with foliations in infinite dimensional manifolds.
In the case of loop groups they are also related with a structure of
manifolds with corners (see the reasons in the introduction).
They are defined with the help of quadrants.
\par {\bf 2.1.1. Remark.} An atlas $At(M)=\{ (U_j,\phi _j): j \} $ 
of a manifold $M$ on a Banach space $X$ over $\bf R$
is called uniform, if its charts satisfy 
the following conditions: \\
$(U1)$ for each $x\in G$ there exist
neighbourhoods $U_x^2\subset U_x^1\subset U_j$
such that for each $y\in U_x^2$ there is the inclusion
$U_x^2\subset U_y^1$; \\
$(U2)$ the image $\phi _j(U_x^2)\subset X$
contains a ball of the fixed positive radius
$\phi _j(U_x^2)\supset B(X,0,r):=\{ y: y\in X,
\| y\| \le r \} ;$ \\
$(U3)$ for each pair of intersecting charts
$(U_1,\phi _1)$ and $(U_2,\phi _2)$ connecting mappings
$F_{\phi _2,\phi _1}=\phi _2\circ \phi _1^{-1}$
are such that $\sup_x \| F'_{\phi _2,\phi _1}(x) \|
\le  C$ and $\sup_x \| F'_{\phi _1,\phi _2}(x) \|
\le  C$, where $C=const >0$ does not depend on  
$\phi _1$ and $\phi _2$. For the diffeomorphism group
$Diff^t_{\beta ,\gamma }(M)$ and loop groups $(L^MN)_{\xi }$
we also suppose that manifolds satisfy conditions
of \cite{ludan,lulgcm,lurim1,lurim2} such that these groups are separable,
but here let $M$ and $N$ may be with a boundary,
where \\
$(N1)$ $N$ is of class not less, than (strongly) $C^{\infty }$ and such that
$sup_{x\in S_{j,l}} \| F^{(n)}_{\psi _j,\psi _l}(x) \| \le C_n$ for each 
$0\le n\in \bf Z$, when $V_{j,l}\ne \emptyset $, $C_n>0$ are constants,
$At(N):= \{ (V_j,\psi _j): j \} $ denotes an atlas of $N$, $V_{j,l}:=
V_j\cap V_l$ are intersections of charts, $S_{j,l}:=\psi _l(V_{j,l})$,
$\bigcup_jV_j=N$. 
\par Conditions $(U1-U3,N1)$ are supposed to be satisfied for the 
manifold $N$ for loop groups, as well as for the manifold $M$ for 
diffeomorphism groups. Certainly, 
the classes of smoothness of manifolds are supposed to be not less 
than that of groups.
\par {\bf 2.1.2.1. Definition.} 
~ A canonical closed subset $Q$ of
$X=\bf R^n$ or of the standard separable
Hilbert space $X=l_2({\bf R})$ over $\bf R$
is called a quadrant if it can be given by $Q:=\{ x\in X:
q_j(x)\ge 0 \} $, where $(q_j: j\in \Lambda _Q)$ 
are linearly independent
elements of the topologically adjoint space $X^*$.  Here 
$\Lambda _Q\subset \bf N$ 
(with $card (\Lambda _Q)=k\le n$ when $X=\bf R^n$)
and $k$ is called the index of $Q$.  If $x\in Q$ and exactly $j$ of the $q_i$'s
satisfy $q_i(x)=0$ then $x$ is called a corner of index $j$.  
Since the unitary space $X=\bf C^n$ or the separable Hilbert space
$l_2({\bf C})$ over $\bf C$ as considered over the field $\bf R$
is isomorphic with $X_{\bf R}:=\bf R^{2n}$ or $l_2({\bf R})$ respectively, 
then the above definition also
describes quadrants in $\bf C^n$ and $l_2({\bf C})$ in such sense. 
In the latter case we also consider 
generalized quadrants
as canonical closed subsets which can be given by
$Q:=\{ x\in X_{\bf R}:$ $q_j(x+a_j)\ge 0, a_j\in X_{\bf R},
j\in \Lambda _Q \} ,$ where $\Lambda _Q\subset \bf N$
($card(\Lambda _Q)=k\in \bf N$ when $dim_{\bf R}X_{\bf R}<\infty $).
\par {\bf 2.1.2.2. Notation.} If for each open subset $U\subset Q\subset X$ 
a function $f:  Q\to Y$ for Banach spaces $X$ and $Y$ over $\bf R$
has continuous Frech\'et differentials $D^{\alpha }f|_U$ on $U$ 
with $\sup_{x\in U} \| D^{\alpha }f(x) \|_{L(X^{\alpha },Y)} <\infty $ 
for each $0\le \alpha
\le r$ for an integer $0\le r$ or $r=\infty $, 
then $f$ belongs to the class of
smoothness $C^r(Q,Y)$, where $0\le r\le \infty $,
$L(X^k,Y)$ denotes the Banach space of bounded $k$-linear
operators from $X$ into $Y$.
\par {\bf 2.1.2.3. Definition.} 
A differentiable mapping $f:  U\to U'$ is called a 
diffeomorphism if 
\par $(i)$ $f$ is bijective and there exist continuous
$f'$ and $(f^{-1})'$, where $U$ and $U'$ 
are interiors of quadrants $Q$
and $Q'$ in $X$. 
\par In the complex case we consider bounded generalized 
quadrants $Q$ and $Q'$ in $\bf C^n$ or $l_2({\bf C})$
such that they are domains with piecewise 
$C^{\infty }$-boundaries and we impose additional conditions on
the diffeomorphism $f$:
\par $(ii)$ ${\bar \partial }f=0$ on $U$, 
\par $(iii)$ $f$ and all its strong (Frech\'et) differentials (as multilinear
operators) are bounded on $U$, 
where $\partial f$ and ${\bar \partial }f$
are differential $(1,0)$ and $(0,1)$ forms respectively,
$d=\partial +{\bar \partial }$ is an exterior derivative.
In particular for 
$z=(z^1,...,z^n)\in \bf C^n$, $z^j\in \bf C$, $z^j=x^{2j-1}+ix^{2j}$
and $x^{2j-1}, x^{2j}\in \bf R$ for each $j=1,...,n,$
$i=(-1)^{1/2}$, there are expressions:
$\partial f:=\sum_{j=1}^n(\partial f/\partial z^j)dz^j$,
${\bar \partial }f:=\sum_{j=1}^n(\partial f/\partial {\bar z}^j)d{\bar z}^j$.
In the infinite dimensional case there are equations:
$(\partial f)(e_j)=\partial f/\partial z^j$
and $({\bar \partial }f)(e_j)=\partial f/\partial {\bar z}^j$, 
where $\{ e_j: {j\in \bf N} \} $ is the standard orthonormal base
in $l_2({\bf C})$, $\partial f/\partial z^j=(\partial f/\partial x^{2j-1}
-i\partial f/\partial x^{2j})/2$,
$\partial f/\partial {\bar z}^j=(\partial f/\partial x^{2j-1}
+i\partial f/\partial x^{2j})/2$.
\par Cauchy-Riemann Condition $(ii)$ means that $f$ on $U$ is the holomorphic
mapping.
\par {\bf 2.1.2.4. Definition and notation.}  
A complex manifold $M$ with corners is
defined in the usual way:  it is a metric separable space
modelled on $X=\bf C^n$ or $X=l_2({\bf C})$ 
and is supposed to be of class $C^{\infty }$.  Charts on $M$ are
denoted $(U_l, u_l, Q_l)$, that is $u_l:  U_l\to u_l(U_l) \subset Q_l$ are
$C^{\infty }$-diffeomorphisms, $U_l$ are open in $M$, $u_l\circ
{u_j}^{-1}$ are biholomorphic from
domains $u_j(U_l\cap U_j)\ne \emptyset $ onto $u_l(U_l\cap U_j)$ 
(that is $u_j\circ u_l^{-1}$ and 
$u_l\circ u_j^{-1}$ are holomorphic and bijective)
and $u_l\circ u_j^{-1}$ 
satisfy conditions $(i-iii)$ from \S 2.1.2.3, $\bigcup_jU_j=M$.
\par A point $x\in M$ is called a corner of index $j$
if there exists a chart $(U,u,Q)$ of $M$ with $x\in U$ and $u(x)$ is of index
$ind_M(x)=j$ in $u(U)\subset Q$.  The set of all corners of index $j\ge 1$ is
called the border $\partial M$ of $M$, $x$ is called an inner point of $M$ if
$ind_M(x)=0$, so $\partial M=\bigcup_{j\ge 1}\partial ^jM$, where
$\partial ^jM:=\{ x\in M:  ind_M(x)=j \} $.  
\par For the real manifold with corners on the connecting mappings
$u_l\circ u_j^{-1}\in C^{\infty }$ 
of real charts is imposed only Condition $2.1.2.3(i)$.
\par {\bf 2.1.2.5. Definition of a submanifold with corners.} 
A subset $Y\subset M$ is called a 
submanifold with corners of $M$ if for each
$y\in Y$ there exists a chart $(U,u,Q)$ of $M$ 
centered at $y$ (that is $u(y)=0$
) and there exists a quadrant $Q'\subset {\bf C^k}$ or in $l_2({\bf C})$
such that
$Q'\subset Q$ and $u(Y\cap U)=u(U)\cap Q'$.  A submanifold with corners $Y$ of
$M$ is called neat, if the index in $Y$ of each $y\in Y$ coincides with its
index in $M$. 
\par  Analogously for real manifolds with corners for $\bf R^k$ and
$\bf R^n$ or $l_2({\bf R})$ instead of $\bf C^k$ and $\bf C^n$
or $l_2({\bf C})$.
\par {\bf 2.1.2.6. Term a complex manifold.} 
Henceforth, the term a complex manifold
$N$ modelled on $X=\bf C^n$ or $X=l_2({\bf C})$ means a metric separable
space supplied with an atlas $\{ (U_j,\phi _j): j\in \Lambda _N \} $
such that:
\par $(i)$ $U_j$ is an open subset of $N$ for each $j\in \Lambda _N$
and $\bigcup_{j\in \Lambda _N}U_j=N$, where $\Lambda _N\subset \bf N$;
\par $(ii)$ $\phi _j: U_j\to \phi _j(U_j)\subset X$ are
$C^{\infty }$-diffeomorphisms, where $\phi _j(U_j)$ are 
$C^{\infty }$-domains in $X$;
\par $(iii)$ $\phi _j\circ \phi _m^{-1}$ is a biholomorphic mapping
from $\phi _m(U_m\cap U_j)$ onto
$\phi _j(U_m\cap U_j)$ while $U_m\cap U_j\ne \emptyset $. 
When $X=l_2({\bf C})$ it is supposed, that
$\phi _j\circ \phi _m^{-1}$ are Frech\'et (strongly)
$C^{\infty }$-differentiable.
\par {\bf 2.1.3.1. Remark.}  Let $X$ be either the standard separable
Hilbert space $l_2=l_2({\bf C})$ over the field $\bf C$ of complex numbers
or $X=\bf C^n$.  Let $t\in \bf
N_o$ $:={\bf N}\cup \{ 0\}$, ${\bf N}:=\{ 1,2,3,... \}$ 
and $W$ be a domain with a continuous piecewise
$C^{\infty }$-boundary $\partial W$
in $\bf R^{2m}$,
$m\in \bf N$, that is $W$ is a $C^{\infty }$-manifold with corners
and it is a canonical closed
subset of $\bf C^m$, $cl(Int(W))=W$, where $cl(V)$ denotes the closure of $V$,
$Int(V)$ denotes the interior of $V$ in the corresponding 
topological space.  As
usually $H^t(W,X)$ denotes the Sobolev space of functions 
$f:  W\to X$ for which
there exists a finite norm
\par $\| f\|_{H^t(W,X)}:=(\sum_{|\alpha |\le t}{\|
D^{\alpha }f\|^2}_{L^2 (W,X)})^{1/2}<\infty $, \\
where $f(x)=(f^j(x):  j\in {\bf
N})$, $f(x)\in l_2$, $f^j(x)\in \bf C$, $x\in W$, 
\par ${\| f\|^2}_{L^2
(W,X)}:=\int_W {\| f(x)\|^2}_X\lambda (dx)$, $\lambda $ is the Lebesgue measure
on $\bf R^{2m}$, $\| z\|_{l_2}:= (\sum_{j=1}^{\infty }|z^j|^2)^{1/2}$, $z=(z^j:
j\in {\bf N})\in l_2$, $z^j\in \bf C$. Then $H^{\infty }(W,X):=
\bigcap_{t\in \bf N}H^t(W,X)$ is the uniform space with the uniformity
given by the family of norms $\{ \| f \|_{H^t(W,X)}: t\in {\bf N} \}$.
\par {\bf 2.1.3.2. Sobolev spaces for manifolds.}  Let now $M$ be a
compact Riemann or complex
$C^{\infty }$-manifold with corners with 
a finite atlas $At(M):=\{ (U_i, \phi _i,
Q_i); i\in \Lambda _M \} $, where $U_i$ are open in $M$, 
$\phi _i:  U_i\to \phi
_i(U_i)\subset Q_i\subset \bf R^m$ (or it is a subset in $\bf C^m$)
are diffeomorphisms
(in addition holomorphic respectively as in \S 2.1.2.3), $(U_i, \phi _i)$ are
charts, $i\in \Lambda _M \subset \bf N$.  
\par Let also $N$ be a
separable complex metrizable manifold with corners 
modelled either on $X=\bf C^n$
or on $X=l_2({\bf C})$ respectively. Let 
$(V_i, \psi _i, S_i)$ be
charts of an atlas $At(N):=\{ (V_i, \psi _i, S_i):  i\in \Lambda _N \} $
such that
$\Lambda _N\subset \bf N$ and $\psi _i:  V_i\to \psi _i(V_i)\subset 
S_i\subset X$ are
diffeomorphisms, $V_i$ are open in $N$, $\bigcup_{i\in \Lambda _N}V_i=N$.
We denote by $H^t(M,N)$ the Sobolev
space of functions $f:  M\to N$ for which $f_{i,j}\in H^t(W_{i,j},X)$ for each
$j\in \Lambda _M$ and $i\in \Lambda _N$ for a domain $W_{i,j}\ne \emptyset $ of
$f_{i,j}$, where $f_{i,j}:=\psi _i\circ f\circ {\phi _j}^{-1}$, and
$W_{i,j}=\phi
_j(U_j\cap f^{-1}(V_i))$ are canonical closed subsets of $\bf R^m$
(or $\bf C^m$ respectively).  The
uniformity in $H^t(M,N)$ is given by the following base $\{ (f,g)\in
(H^t(M,N))^2:  \sum_{i\in \Lambda _N, j\in \Lambda _M} {\| f_{i,j}-
g_{i,j}\|^2}_{H^t(W_{i,j}, X)}<\epsilon \}$, where $\epsilon >0$, $W_{i,j}$
are domains of $(f_{i,j}- g_{i,j}).$ For $t=\infty $ as usually
$H^{\infty }(M,N):=\bigcap_{t\in \bf N} H^t(M,N)$.
\par {\bf 2.1.3.3. A uniform space of piecewise holomorphic mappings.}  
For two complex manifolds $M$ and $N$ 
with corners let 
${\sf O}_{\Upsilon }(M,N)$ denotes a space
of continuous mappings $f: M\to N$ such that
for each $f$ there exists a partition $Z_f$ of $M$ with the help of
a real $C^{\infty }$-submanifold ${M'}_f$, 
which may be with corners, such that its codimension 
over $\bf R$ in $M$ is $codim_{\bf R}{M'}_f=1$
and $M\setminus {M'}_f$ is a disjoint union of open complex submanifolds
$M_{j,f}$ possibly with corners 
with $j=1,2,...$ such that each restriction $f|_{M_{j,f}}$
is holomorphic with all its derivatives bounded on $M_{j,f}$.
For a given partition $Z$ (instead of $Z_f$) and the corresponding
$M'$ the latter subspace of continuous piecewise holomorphic mappings
$f: M\to N$ is denoted by ${\sf O}(M,N;Z)$.
The family $\{ Z \} $ of all such partitions is denoted
$\Upsilon $.
That is ${\sf O}_{\Upsilon } (M,N)=str-ind_{\Upsilon }{\sf O}(M,N;Z)$.
Let also ${\sf O}(M,N)$ denotes the space
of holomorphic mappings $f: M\to N$,
$Diff^{\infty }(M)$ denotes a group of 
$C^{\infty }$-diffeomorphisms of $M$ and $Diff^{\sf O}_{\Upsilon }(M):=
Hom(M)\cap {\sf O}_{\Upsilon }(M,M)$, where $Hom(M)$ is a group
of homeomorphisms.
\par Let $A$ and $B$ be two complex manifolds with corners
such that $B$ is a submanifold of $A$. Then $B$ is called
a strong $C^r([0,1]\times A,A)$-retract 
(or $C^r([0,1],{\sf O}_{\Upsilon }(A,A))$-retract) of  $A$ if there exists
a mapping $F: [0,1]\times A\to A$ such that $F(0,z)=z$ for each
$z\in A$ and $F(1,A)=B$ and $F(x,A)\supset B$ for each $x\in
[0,1]:=\{ y: 0\le y \le 1, y\in {\bf R} \} $,
$F(x,z)=z$ for each $z\in B$ and $x\in [0,1]$,
where $F\in C^r([0,1]\times A,A)$ or $F\in C^r([0,1],{\sf O}_{\Upsilon }
(A,A))$ respectively, $r\in [0,\infty )$, $F=F(x,z)$, $x\in [0,1]$,
$z\in A$. Such $F$ is called the retraction. In the case of $B=\{ a_0 \} $,
$a_0\in A$ we say that $A$ is $C^r([0,1]\times A,A)$-contractible
(or $C^r([0,1], {\sf O}_{\Upsilon }(A,A))$-contractible correspondingly).
Two maps $f: A\to E$ and $h: A\to E$ are called 
$C^r([0,1]\times A,E)$-homotopic (or 
$C^r([0,1],{\sf O}_{\Upsilon }(A,E))$-homotopic ) if there
exists $F\in C^r([0,1]\times A,E)$ (or $F\in C^r([0,1],{\sf O}_{\Upsilon }
(A,E))$ respectively) such that $F(0,z)=f(z)$ and $F(1,z)=h(z)$
for each $z\in A$, where $E$ is also a complex manifold.
Such $F$ is called the homotopy.
\par Let $M$ be a complex manifold with corners
satisfying the following conditions:
\par $(i)$ it is compact;
\par $(ii)$ $M$ is a union of two closed complex submanifolds $A_1$ and $A_2$
with corners,
which are canonical closed subsets in $M$
with $A_1\cap A_2=\partial A_1\cap \partial A_2=:A_3$ 
and a codimension over $\bf R$ of $A_3$ in $M$ is $codim_{\bf R}A_3=1$;
\par $(iii)$ a marked point $s_0$ is in $A_3$;
\par $(iv)$ $A_1$  and $A_2$ are 
$C^0([0,1],{\sf O}_{\Upsilon }(A_j,A_j))$-contractible into a marked point
$s_0\in A_3$ by mappings $F_j(x,z),$
where either $j=1$ or $j=2$.
There can be considered more general condition
of $C^0([0,1],{\sf O}_{\Upsilon }(A_j,A_j))$-contractibility
of $A_j$ on $X_0\cap A_j$, where $X_0$ is a closed subset in $M$,
$j=1$ or $j=2$, $s_0\in X_0$.
\par  We consider all finite partitions $Z:=\{ M_k:  k\in \Xi _Z\}$ of $M$
such that $M_k$ are complex submanifolds (of $M$), which may be
with corners and
$\bigcup_{k=1}^sM_k=M$, $\Xi _Z=\{ 1,2,...,s \}$, $s\in \bf N$ 
depends on $Z$, $M_k$ are
canonical closed subsets of $M$.  We denote by $\tilde diam(Z):=\sup_k(diam
(M_k))$ the diameter of the partition $Z$, where $diam (A)=\sup_{x, y\in A}
|x-y|_{\bf C^n}$ is a diameter of a subset $A$ in $\bf C^n$, since
each finite dimensional manifold $M$ can be embedded into $\bf C^n$
with the corresponding $n\in \bf N$.  
We suppose
also that $M_i\cap M_j\subset M'$ and $\partial M_j\subset M'$ for each $i\ne
j$, where $M'$ is a closed $C^{\infty }$-submanifold (which may be
with corners) in $M$ with
the codimension $codim_{\bf R}(M')=1$ of $M'$ in $M$, $M'=\bigcup_{j\in \Gamma
_Z}{M'}_j$, ${M'}_j$ are $C^{\infty }$-submanifolds of $M$, 
$\Gamma _Z$ is a finite
subset of $\bf N$. 
\par  We denote by $H^t(M,N;Z)$ a space of continuous functions
$f:  M\to N$ such that $f|_{(M\setminus M')}\in H^t(M\setminus M',N)$ and
$f|_{[Int(M_i)\cup (M_i\cap {M'}_j)]}\in H^t(Int(M_i)\cup (M_i\cap {M'}_j),N)$,
when $\partial M_i\cap {M'}_j\ne \emptyset $, $h^Z_{Z'}:  H^t(M,N;Z)\to
H^t(M,N;Z')$ are embeddings for each $Z\le Z'$ in $\Upsilon $.
\par The ordering
$Z\le Z'$ means that each submanifold $M_i^{Z'}$
from a partition $Z'$ either belongs to the family
$(M_j:  j=1,...,k)=(M_j^Z: j=1,...,k)$ for $Z$ or 
there exists $j$ such that $M_i^{Z'}\subset M_j^Z$ and
$M_j^Z$ is a finite union of $M_l^{Z'}$ for which $M_l^{Z'}\subset
M_j^Z$. Moreover, these $M_l^{Z'}$ are submanifolds (may be
with corners) in $M_j^Z$. 
\par Then we consider the following uniform space
$H^t_p(M,N)$ that is the strict inductive limit $str-ind \{ H^t(M,N;Z);
h^{Z'}_{Z}; \Upsilon \} $ (the index $p$ reminds about the procedure of
partitions), where $\Upsilon $ is the directed family of all such $Z$, 
for which $\lim_{\Upsilon }\tilde diam(Z)=0$.  
\par {\bf 2.1.4. Notes and definitions of loop 
monoids and loop groups.}  Let now $s_0$ be the marked
point in $M$ such that $s_0\in A_3$ (see \S 2.1.3.3) 
and $y_0$ be a marked point in the manifold $N$.
\par $(i).$ Suppose that $M$ and $N$ are connected.
\par Let
$H^t_p(M,s_0;N,y_0):=\{ f\in H^t(M,N)| f(s_0)=y_0 \} $ denotes the closed
subspace of $H^t(M,N)$ and $\omega _0$ be its element such that $\omega
_0(M)= \{ y_0 \} $, where $\infty \ge t\ge m+1$, $2m=dim_{\bf R}M$
such that $H^t\subset C^0$ due to the Sobolev embedding theorem.  
The following subspace $ \{ f: f\in H^{\infty }_p(M,s_0;N,y_0),\mbox{ }
\mbox{ } {\bar \partial }f=0 \} $ is isomorphic with
${\sf O}_{\Upsilon }(M,s_0;N,y_0)$, since $f|_{(M\setminus M')}
\in H^{\infty }(M\setminus M',N)=
C^{\infty }(M\setminus M',N)$ and ${\bar \partial }f=0$.
\par Let as usually $A\vee B:=A\times \{
b_0\}\cup \{ a_0\}\times B\subset A\times B$ 
be the wedge sum of 
pointed spaces $(A,a_0)$ and $(B,b_0)$, where $A$
and $B$ are topological spaces with marked points $a_0\in A$ and $b_0\in B$.
Then the wedge combination $g\vee f$ 
of two elements $f, g\in H^t_p(M,s_0;N,y_0)$ is
defined on the domain $M\vee M$.
\par The spaces ${\sf O}_{\Upsilon }(J,A_3;N,y_0):=\{ f\in {\sf O}_{\Upsilon }
(J,N): f(A_3)=\{ y_0 \} \} $ 
have the manifold structure and have embeddings into
${\sf O}_{\Upsilon }(M,s_0;N,y_0)$ due to Condition 
$2.1.3.3(ii)$, where either
$J=A_1$ or $J=A_2$.
This induces the following embedding
$\chi ^*: {\sf O}_{\Upsilon }(M\vee M,s_0\times s_0;N,y_0)\hookrightarrow
{\sf O}_{\Upsilon }(M,s_0;N,y_0)$.
Considering $H^t_p(M,X_0;N,y_0)=\{ f\in H^t(M,N): f(X_0)=\{ y_0 \} \} $
and ${\sf O}_{\Upsilon }(J,A_3\cup X_0;N,y_0)$ we get
the embedding 
$\chi ^*: {\sf O}_{\Upsilon }(M\vee M,X_0\times X_0;N,y_0)\hookrightarrow
{\sf O}_{\Upsilon }(M,X_0;N,y_0)$.
Therefore $g\circ f:=\chi ^*(f\vee g)$ is the composition in
${\sf O}_{\Upsilon }(M,s_0;N,y_0)$.  
\par The space $C^{\infty }(M,N)$ is dense in $C^0(M,N)$
and there is the inclusion ${\sf O}_{\Upsilon }(M,N)\subset
H^{\infty }_p(M,N)$. Let $M_{\bf R}$ be the Riemann manifold generated by
$M$ considered over $\bf R$. Then $Diff^{\infty }_{s_0}(M_{\bf R})$
is a group of $C^{\infty }$-diffeomorphisms $\eta $ of $M_{\bf R}$ 
preserving the marked point $s_0$, that is $\eta (s_0)=s_0$.
There exists the following equivalence relation
$R_{\sf O}$ in ${\sf O}_{\Upsilon }(M,s_0;N,y_0)$:
$fR_{\sf O}h$ if and only if there exist nets
$\eta _n\in Diff^{\infty }_{s_0}(M_{\bf R})$, 
also $f_n$ and $h_n\in H^{\infty }_p(M,s_0;N,y_0)$ with $\lim_{n}f_n=f$ 
and $\lim_{n}h_n=h$ such that $f_n(x)=h_n(\eta _n(x))$ for each
$x\in M$ and $n\in \omega $, where $\omega $ is a directed set,
$f, h \in {\sf O}_{\Upsilon }(M,s_0;N,y_0)$ and converegence is considered
in $H^{\infty }_p(M,s_0;N,y_0)$. 
In general case we consider $Diff^{\infty }_{X_0}(M_{\bf R}):=
\{ f\in Diff^{\infty }(M_{\bf R}): f(X_0)=X_0 \} $ and elements
$f, h$ in ${\sf O}_{\Upsilon }(M,X_0;N,y_0)$ and convergence
in $H^{\infty }(M,X_0;N,y_0)$ we get the equivalence relation
$R_{\sf O}$ in ${\sf O}_{\Upsilon }(M,X_0;N,y_0)$.
\par The quotient space ${\sf O}_{\Upsilon }(M,X_0;N,y_0)/R_{\sf O}=:
(S^MN)_{\sf O}$ is called the loop monoid. 
It has a structure of topological Abelian monoid
with the cancellation property (see \cite{ludan,lulgcm}).
Applying the A. Grothendieck procedure (see below)
to $(S^MN)_{\sf O}$ we get a loop group $(L^MN)_{\sf O}$.
For the spaces $H^t_p(M,s_0;N,y_0)$ the corresponding equivalence 
relations are denoted $R_{t,H}$, the loop monoids are denoted 
by $(S^M_{\bf R}N)_{t,H}$, the loop groups are denoted
by $(L^M_{\bf R}N)_{t,H}$. When real  manifolds are considered 
we omit the index $\bf R$.
\par For a commutative 
monoid with the cancellation property
$(S^MN)_{\sf O}$ there exists a 
commutative group $(L^MN)_{\sf O}$
equal to the Grothendieck group. This group algebraically is the quotient group
$F/\sf B$, where $F$ is a free Abelian group generated by 
$(S^MN)_{\sf O}$ and $\sf B$ is a subgroup of $F$ generated by
elements $[f+g]-[f]-[g]$, $f$ and $g\in (S^MN)_{\sf O}$,
$[f]$ denotes an element of $F$ corresponding to $f$. The natural mapping 
$\gamma : (S^MN)_{\sf O}\to (L^MN)_{\sf O}$ is injective.
We supply $F$ with a topology inherited from
the Tychonoff product topology of $(S^MN)_{\sf O}^{\bf Z}$,
where each element $z$ of $F$ is $z=\sum_fn_{f,z}[f]$, 
$n_{f,z}\in \bf Z$ for each $f\in (S^MN)_{\sf O}$,
$\sum_f|n_{f,z}|<\infty $. In particular $[nf]-n[f]\in \sf B$,
hence $(L^MN)_{\sf O}$ is the complete topological group and
$\gamma $ is the topological embedding
such that $\gamma (f+g)=\gamma (f)+ \gamma (g)$
for each $f, g \in (S^MN)_{\sf O}$, $\gamma (e)=e$,
since $(z+B)\in \gamma (S^MN)_{\sf O}$, when $n_{f,z}\ge 0$
for each $f$, so in general $z=z^{+}-z^{-}$, where
$(z^{+}+B)$ and $(z^{-}+B)\in \gamma (S^MN)_{\sf O}$.
\par {\bf 2.1.5.1. Notes and Definitions.} In view of \S I.5 \cite{kob}
a complex manifold $M$ considered over $\bf R$
admits a (positive definite) Riemann metric $g$, since $M$ is paracompact
(see \S \S 1.3 and 1.5 \cite{kob}). Due to Theorem IV.2.2
\cite{kob} there exists the Levi-Civit\`a connection (with vanishing torsion)
of $M_{\bf R}$. 
For the orientable manifold $M$ suppose $\nu $ is a measure on 
$M$ corresponding to the Riemann volume element $w$ ( $m$-form )
$\nu (dx)=w(dx)/w(M)$. The Riemann volume element $w$ is non-degenerate
and non-negative, since $M$ is orientable.
For the nonorientable $M$
consider $\tilde M$ its double covering orientable manifold 
and the quotient mapping $\theta _M: \tilde M\to M$,
then the Riemann volume element $w$ on $\tilde M$
produces the following measure $\nu (S):=w(\theta _M^{-1}(S))/
(2w(\tilde M))$ for each Borel subset $S$ in $M$.
\par The Christoffel
symbols $\Gamma ^k_{i,j}$ of the Levi-Civit\`a derivation
(see \S 1.8.12 \cite{kling}) are of class $C^{\infty }$ for $M$.
Then the
equivalent uniformity in $H^t(M,N)$ for $0\le t<\infty $
is given by the following base $\{
(f,g)\in (H^t(M,N))^2:  {\| (\psi _j\circ f -\psi _j\circ g)\|"}_{
H^t(M,X)}<\epsilon $, 
where $D^{\alpha }=\partial ^{|\alpha |}/\partial (x
^1)^{\alpha ^1} ...\partial (x ^{2m})^{\alpha ^{2m}}$, $\epsilon >0$, 
${{\| (\psi
_j\circ f -\psi _j\circ g)\|"}^2}_{ H^t(M,X)}:= \sum_{|\alpha |\le t}
\int_{M}|D^{\alpha } (\psi _j\circ f(x ) 
-\psi _j\circ g(x ))|^2\nu
(dx) \} $, $j\in \Lambda _N$, $X$ is the Hilbert space  over $\bf C$
either $\bf C^n$ or $l_2({\bf C})$, 
$x$ are local normal coordinates in $M_{\bf R}$.  
We consider submanifolds $M_{i,k}$ and ${M'}_{j,k}$
for each partition $Z_k$ as in \S 2.1.3.3 
(with $Z_k$ instead of $Z$), $i\in \Xi
_{Z_k}$, $j\in \Gamma _{Z_k}$, where 
$\Xi _{Z_k}$ and $\Gamma _{Z_k}$ are finite
subsets of $\bf N$. 
We supply $H^{\gamma }(M,X;Z_k)$ with the following metric 
$\rho _{k, \gamma }(y):=[\sum_{i\in
\Xi }$ $ {{\| y|_{M_{i,k}}\|"}^2}_{\gamma ,i,k}]^{1/2}$ for $y\in H^{\gamma
}(M,X;Z_k)$ and $\rho _{k, \gamma }(y)=+\infty $ in the 
contrary case, where $\Xi =\Xi
_{Z_k}$, $\infty >t\ge \gamma \in \bf N$,
$\gamma \ge m+1$, ${\|
y|_{M_{i,k}}\|"}_{\gamma ,i,k}$ is given analogously to ${\| y\|"}_{H^{\gamma
}(M,X)}$, but with $\int_{M_{i,k}}$ instead of $\int_{M}$.  
\par Let $Z^{\gamma
}(M,X)$ be the completion of $str-ind \{ H^{\gamma }(M,X;Z_j); h^{Z_i}_{Z_j};
{\bf N} \}=:Q$ relative to the following norm 
${\| y\| '}_{\gamma }:= \inf_k\rho
_{k, \gamma }(y)$, as usually
$Z^{\infty }(M,X)=\bigcap_{\gamma \in \bf N}Z^{\gamma }(M,X)$.
Let ${\bar Y}^{\infty }(M,X):= \{ f: f\in Z^{\infty }(M,X),$
${\bar \partial }f_j|_{M_{j,k}}=0$ $\mbox{for each }k \} ,$
where $f\in Z^{\infty }(M,X)$ imples $f=\sum_jf_j$ with
$f_j\in H^{\infty }(M,X;Z_j)$ for each $j\in \bf N$.
\par For a domain $W$ in $\bf C^m$, which is a complex manifold with corners,
let $Y^{\Upsilon ,a}(W,X)$ (and $Z^{\Upsilon ,a}(W,X)$)
be a subspace of those
$f\in {\bar Y}^{\infty }(W,X)$ (or $f\in Z^{\infty }(W,X)$ respectively)
for which 
$$\| f\|_{\Upsilon ,a}:=(\sum_{j=0}^{\infty }
({\| f\|^{*}}_j)^2/[(j!)^{a_1}j^{a_2}])^{1/2}<\infty ,$$ 
where ${({\| f\|^{*}}_j)^2}:=
({\| f\| '}_j)^2-
({\| f\| '}_{j-1})^2$ for $j\ge 1$ and ${{\| f\|^*}_0}={\| f\|'}_0$, 
$a=(a_1,a_2)$, $a_1$ and $a_2\in \bf R$, $a<a'$ if either
$a_1<{a'}_1$ or $a_1={a'}_1$ and $a_2<{a'}_2$.
\par Using the atlases $At(M)$ and $At(N)$ 
for $M$ and $N$ of class of smoothness $Y^{\Upsilon ,b}\cap C^{\infty }$ 
with $a\ge b$
we get the uniform space $Y^{\Upsilon ,a}(M,X_0;N,y_0)$ (and also
$Z^{\Upsilon ,a}(M,X_0;N,y_0)$)
of mappings $f:  M\to N$ with $f(X_0)=y_0$ such
that $\psi _j\circ f\in Y^{\Upsilon ,a}(M,X)$ 
(or $\psi _j\circ f\in Z^{\Upsilon ,a}(M,X)$ respectively) for each $j$,
where $\sum_{p\in \Lambda _M, j\in \Lambda _N}\| f_{p,j}-(w_0)_{p,j}
\|^2_{Y^{\Upsilon ,a}(W_{p,j},X)}<\infty $ for each 
$f\in Y^{\Upsilon ,a}(M,X_0;N,y_0)$ is satisfied with
$w_0(M)=\{ y_0 \} $, since $M$ is compact.
Substituting $w_0$ on a fixed mapping
$\theta : M\to N$ we get the uniform space
$Y^{\Upsilon ,a,\theta }(M,N)$.
To each equivalence class $\{ g:  gR_{\sf O}f \} 
=:<f>_{\sf O}$ there corresponds an
equivalence class $<f>_{\Upsilon ,a}:= cl(<f>_{\sf O}\cap Y^{\Upsilon ,a}
(M,X_0;N,y_0))$
(or $<f>^{\bf R}_{\Upsilon ,a}:= cl(<f>_{\infty ,H}\cap Z^{\Upsilon ,a}
(M,X_0;N,y_0))$), where the closure is taken
in $Y^{\Upsilon ,a}(M,X_0;N,y_0)$ 
(or $Z^{\Upsilon ,a}(M,X_0;N,y_0)$ respectively). 
This generates equivalence relations $R_{\Upsilon ,a}$
and $R^{\bf R}_{\Upsilon ,a}$ respectively.
We denote the quotient spaces
$Y^{\Upsilon ,a}(M,X_0;N,y_0)/R_{\Upsilon ,a}$ 
and $Z^{\Upsilon ,a}(M,X_0;N,y_0)/R^{\bf R}_{\Upsilon ,a}$ 
by $(S^MN)_{\Upsilon ,a}$ and $(S^M_{\bf R}N)_{\Upsilon ,a}$
correspondingly. Using the A. Grothendieck construction we get
the loop groups $(L^MN)_{\Upsilon ,a}$ and $(L^M_{\bf R}N)_{\Upsilon ,a}$
respectively.
\par {\bf 2.1.5.2. Gevrey-Sobolev classes of smoothness
of loop monoids and loop groups. Notes and definitions.} 
Let $M$ be an infinite dimensional complex 
$Y^{\xi '}$-manifold
with corners modelled on $l_2({\bf C})$ such that 
\par $(i)$ there is the sequence 
of the canonically embedded
complex submanifolds $\eta _m^{m+1}:
M_m\hookrightarrow M_{m+1}$ for each $m\in \bf N$
and to $s_{0,m}$ in $M_m$ it corresponds $s_{0,m+1}=
\eta _m^{m+1}(s_{0,m})$ in $M_{m+1}$, $dim_{\bf C}M_m=n(m)$,
$0<n(m)<n(m+1)$ for each $m\in \bf N$, $\bigcup_mM_m$ is dense in $M$;
\par $(ii)$ $M$ and $At(M)$ are foliated, that is, 
\par $(\alpha )$ $u_i\circ u_j^{-1}|_{u_j(U_i\cap U_j)}\to l_2$
are of the form: $u_i\circ u_j^{-1}((z^l:$ $l\in {\bf N}))=
(\alpha _{i,j,m}(z^1,...,z^{n(m)}), \gamma _{i,j,m}(z^l:$ $l>n(m)))$
for each $m$, when $M$ is without a boundary. 
If $\partial M\ne \emptyset $ then 
\par $(\beta )$ for each boundary component $M_0$ of $M$ 
and $U_j\cap M_0\ne \emptyset $ we have $\phi _j: U_j\cap M_0\to
H_{l,Q}$, moreover, $\partial M_m\subset \partial M$ for each $m$,
where $H_{l,Q}:= \{ z\in Q_j:$ $x^{2l-1}\ge 0 \} $,
$Q_j$ is a quadrant in $l_2$ such that $Int_{l_2}H_{l,Q}\ne \emptyset $
(the interior of $H_{l,Q}$ in $l_2$), $z^l=x^{2l-1}+ix^{2l}$,
$x^j\in \bf R$, $z^l\in \bf C$ (see also \S 2.1.2.4);
\par $(iii)$ $M$ is embedded into $l_2$ as a bounded closed subset;
\par $(iv)$ each $M_m$ satisfies conditions $2.1.3.3 (i-iv)$
with $X_{0,m}:=X_0\cap M_m$, where $X_0$ is a closed subset in $M$.
\par Let $W$ be a bounded canonical closed subset
in $l_2({\bf C})$ with a continuous piecewise $C^{\infty }$-boundary and $H_m$ 
an increasing sequence of finite dimensional subspaces
over $\bf C$, $H_m\subset H_{m+1}$ and $dim_{\bf C}H_m=n(m)$ for
each $m\in \bf N$. Then there are spaces 
$P^{\infty }_{\Upsilon ,a}(W,X):=str-ind_mY^{\Upsilon ,a}(W_m,X)$,
where $W_m=W\cap H_m$ and $X$ is a separable Hilbert space over $\bf C$.
\par Let $Y^{\xi }(W,X)$
be the completion of $P^{\infty }_{\Upsilon ,a}(W,X)$
relative to the following norm 
$$\| f\|_{\xi }:=
[\sum_{m=1}^{\infty }{{\| f|_{W_m} \| " }^2}_{Y^{\Upsilon ,a}(W_m,X)}/[
(n(m)!)^{1+c_1}n(m)^{c_2}]]^{1/2},$$
where 
${{\| f|_{W_m} \| " }^2}_{Y^{\Upsilon ,a}(W_m,X)}:=
{\| f|_{W_m} \| ^2}_{Y^{\Upsilon ,a}(W_m,X)}
-{\| f|_{W_{m-1}} \| ^2}_{Y^{\Upsilon ,a}(W_{m-1},X)}$
for each $m>1$ and
${\| f|_{W_1} \| "}_{Y^{\Upsilon ,a}(W_1,X)}:=
\| f|_{W_1} \| _{Y^{\Upsilon ,a}(W_1,X)}$;
$c=(c_1,c_2)$, $c_1$ and $c_2\in \bf R$, $c<c'$ if either 
$c_1<{c'}_1$ or $c_1={c'}_1$ and $c_2<{c'}_2$; $\xi =(\Upsilon ,a,c)$.
Let $M$ and $N$ be the $Y^{\Upsilon ,a',c'}$-manifolds with 
$a'<a$ and $c'<c$.
\par If $N$ is the finite dimensional complex $Y^{\Upsilon ,a'}$-manifold,
then it is also the $Y^{\Upsilon ,a',c'}$-manifold.
There exists the strict inductive limit of loop groups
$(L^{M_m}N)_{\Upsilon ,a}=:L^m$, since
there are natural embeddings $L^m\hookrightarrow L^{m+1}$,
such that each element $f\in Y^{\Upsilon ,a}(M_m,X_{0,m};N,y_0)$ is
considered in $Y^{\Upsilon ,a}(M_{m+1},X_{0,m+1};N,y_0)$
as independent from $(z^{n(m)+1},...,z^{n(m+1)-1})$ in the local normal
coordinates $(z^1,...,z^{n(m+1)})$ of $M_{m+1}$.
We denote it $str-ind_mL^m=:(L^MN)_{\Upsilon ,a}$
and also $str-ind_mQ^m=:Q^{\infty }_{\Upsilon ,a}(N,y_0)$,
\par $str-ind_mY^{\Upsilon ,a}(M_m;N)=:Q^{\infty }_{\Upsilon ,a}(N)$,
where $Q^m:=Y^{\Upsilon ,a}(M_m,X_{0,m};N,y_0)$.
Then with the help of charts of $At(M)$ and $At(N)$ the space
$Y^{\xi }(W,X)$
induces the uniformity $\tau $ in $Q^{\infty }_{\Upsilon ,a}(N,y_0)$
and the completion of it relative to $\tau $ we denote by
$Y^{\xi }(M,X_0;N,y_0)$, where $\xi =(\Upsilon ,a,c)$ and
$\sum_{p\in \Lambda _M, j\in \Lambda _N}\| f_{p,j}-(w_0)_{p,j}
\|^2_{Y^{\xi }(W_{p,j},X)}<\infty $ for each 
$f\in Y^{\xi }(M,X_0;N,y_0)$ is supposed to be satisfied with
$w_0(M)=\{ y_0 \} $, since each $M_m$ is compact.
Substituting $w_0$ on the fixed mapping $\theta : M\to N$ we get the 
uniform space $Y^{\xi ,\theta }(M,N)$.
Therefore,
using classes of equivalent elements from $Q^{\infty }_{\Upsilon ,a}(N,y_0)$ 
and their closures in $Y^{\xi }(M,X_0;
N,y_0)$ as in \S 2.1.5.1 we get
the corresponding loop monoids which are denoted 
$(S^MN)_{\xi }$.  With the help of A. Grothendieck construction we get
loop groups $(L^MN)_{\xi }$. Substituting spaces $Y^{\Upsilon ,a}$
over $\bf C$ 
onto $Z^{\Upsilon ,a}$ over $\bf R$ with respective modifications we get
spaces $Z^{\Upsilon ,a,c}(M,N)$ over $\bf R$,
loop monoids $(S^M_{\bf R}N)_{\xi }$ and groups $(L^M_{\bf R}N)_{\xi }$
for the multi-index $\xi =(\Upsilon ,a,c)$.
\par Let $exp:  {\tilde T}N\to N$ be the exponential mapping, 
where ${\tilde T}N$ is a neighbourhood of $N$ in $TN$
\cite{kling}.
\par The relation between manifolds with corners and usual 
manifolds is given by the following lemma.
\par {\bf 2.1.6. Lemma.} {\it If $M$ is a complex manifold 
modelled on $X=\bf C^n$ or $X=l_2({\bf C})$ 
with an atlas $At(M)=\{ (V_j,\phi _j): j \}$,
then there exists an atlas $At'(M)= \{ (U_k,u_k,Q_k): k \} $
which refines $At(M)$, where $(V_j,\phi _j)$ are usual charts
with diffeomorphisms $\phi _j: V_j\to \phi _j(V_j)$ such that 
$\phi _j(V_j)$ are $C^{\infty }$-domains in $\bf C^n$ and 
$(U_k,u_k,Q_k)$ are charts corresponding to quadrants $Q_k$ in $\bf C^n$
or $l_2({\bf C})$} (see \cite{lulgcm} \S 2.3.1 and \cite{ludan}).
\par Necessary data about structures of loop groups are given in 
Theorems 2.1.7 and 2.1.8.
\par {\bf 2.1.7. Theorems. (1).} {\it The space $(L^MN)_{\xi }=:G$ 
for $\xi =(\Upsilon ,a)$ or $\xi =(\Upsilon ,a,c)$ from \S 2.1.5
is the complete separable Abelian topological group. 
Moreover, $G$ is the dense subgroup in 
$(L^MN)_{\sf O}$ for $\xi =(\Upsilon ,a)$;
$G$ is non-discrete non-locally compact and locally connected.}
\par {\bf (2).} {\it  The space $X^{\xi }(M,N):=T_e(L^MN)_{\xi }$  
is Hilbert for each $1\le m=dim_{\bf C}M\le \infty $.}
\par {\bf (3.)} {\it Let $N$ be a complex Hilbert 
$Y^{\xi '}$-manifold with $a>a'$ and $c>c'$ 
for $\xi '=(\Upsilon ,a')$ or $\xi '=(\Upsilon ,a',c')$ respectively,
then there exists a mapping 
$\tilde E: {\tilde T}(L^MN)_{\xi }\to (L^MN)_{\xi }$ 
defined by $\tilde E_{\eta
}(v)=exp_{\eta (s)}\circ v_{\eta }$ on a neighbourhood $V_{\eta }$ 
of the zero
section in $T_{\eta }(L^MN)_{\xi }$ and it is 
a $C^{\infty }$-mapping for
$Y^{\xi '}$-manifold $N$ by $v$ onto
a neighbourhood $W_{\eta }=W_e\circ \eta $ of $\eta \in (L^MN)_{\xi }$;
$\tilde E$ is the uniform isomorphism of 
uniform spaces $V_{\eta }$ and $W_{\eta
}$, where $s\in M$, $e$ is the unit element in $G$, $v\in V_{\eta },$
$1\le m\le \infty $.}
\par {\bf (4).} {\it $(L^MN)_{\xi }$ is the 
closed proper subgroup in $(L^M_{\bf R}N)_{\xi }$.}
\par The {\bf proof} for the orientable manifolds
follows from \S 2.9 \cite{lulgcm} and \cite{ludan}, 
the case of nonorientable manifolds is analogous due to \S \S 2.1.4
and 2.1.5. The latter case can be deduced also
from Theorems 2.1.8 below and the case of orientable manifolds.
\par {\bf 2.1.8. Theorems.} {\it Suppose that manifolds
$M$ and $N$ together with their covering manifolds 
$\tilde M$ and $\tilde N$ satisfy the conditions imposed above.
\par $(1)$. Let $N$ be the nonorientable manifold and 
$\theta _N: \tilde N\to N$ is the quotient mapping of its
double covering manifold $\tilde N$. 
Then there exists a quotient group homomorphism 
$\tilde \theta _N : (L^M\tilde N)_{\xi }\to (L^MN)_{\xi }$.
\par $(2)$. Let $M$ be the nonorientable manifold,
then the quotient mapping $\theta _M: \tilde M\to M$ induces the 
group embedding $\tilde \theta _M:
(L^MN)_{\xi }\hookrightarrow (L^{\tilde M}N)_{\xi }.$}
\par {\bf Proof.} If $M$ is the nonorientable manifold,
then there exists the homomorphism $h$ of the fundamental group
$\pi _1(M,s_0)$ onto the two-element group $\bf Z_2$.
For connected $M$ the group $\pi _1(M,s_0)$ does not depend 
on the marked point $s_0$ and it is denoted by $\pi _1(M)$.
If $M$ is the connected manifold, 
then it has a universal covering manifold $M^*$ which is 
linearly connected and a fiber bundle with the group
$\pi _1(M)$ and a projection $p: M^*\to M$.
Using the homomorphism $h$ one gets the orientable double covering
$\tilde M$ of $M$ such that $\tilde M$ is connected, if
$M$ is connected (see Proposition 5.9 \cite{kob} and Theorem 78
\cite{pont}). Moreover, for each $x\in M$ there exists
a neighborhood $U$ of $x$ such that 
$\theta _M^{-1}(U)$ is the disjoint union of two diffeomorphic
open subsets $V_1$ and $V_2$ in $\tilde M$, where $\theta _M: \tilde M\to M$
is the quotient mapping, $g: V_1\to V_2$ is a diffeomorphism. 
\par $(1)$.
For each $\tilde f \in Z^{\Upsilon ,a,c}(M,\tilde N)$
there exists $f=\theta _N\circ \tilde f$ in $Z^{\Upsilon ,a,c}(M,N)$.
This induces the quotient mapping $\bar \theta _N :
Z^{\Upsilon ,a,c}(M,\tilde N)\to Z^{\Upsilon ,a,c}(M,N)$,
hence it induces the quotient mapping
$\bar \theta _N :
Z^{\Upsilon ,a,c}(M\vee M,\tilde N)\to Z^{\Upsilon ,a,c}(M\vee M,N)$
such that $\bar \theta _N (f\vee h) =
\bar \theta _N(f)\vee \bar \theta _N(g)$.
Considering the equivalence relation in $Y^{\xi }(M,s_0;\tilde N,y_0)$
and then loop monoids we get the quotient homomorphism
$(S^M\tilde N){\xi }\to (S^MN)_{\xi }$. With the help of A. Grothendieck 
construction it induces the loop groups quotient homomorphism.
\par $(2)$. On the other hand, let $M$ be the nonorientable manifold
then the quotient mapping $\theta _M: \tilde M\to M$ induces the 
locally finite open covering $\{ U_x: x\in M_0 \} $ of $M$,
where $M_0$ is a subset of $M$, such that each $\theta _M^{-1}(U_x)$
is the disjoint union of two open subsets $V_{x,1}$ and $V_{x_2}$ 
in $\tilde M$ and there exists a diffeomorphism $g_x$ of $V_{x,1}$
on $V_{x,2}$ of the same class of smoothness as $M$.
This produces the closed subspace of all $\tilde f\in  Z^{\Upsilon ,a,c}(
\tilde M,N)$ for which 
\par $(i)$ $\tilde f|_{V_{x,1}}(g_x^{-1}(y))=\tilde f|_{V_{x,2}}(y)$
for each $y\in V_{x,2}$ and for each $x\in M_0$, 
where $M_0=M_0^f$ and $\{ U_x=U_x^f: x\in M_0 \} $
may depend on $f$. 
If $s_0$ is a marked point in $M$, then one of the points $\tilde s_0$
of $\theta _M^{-1}(s_0)$ let be the marked point in $\tilde M$.
Then $\theta _M$ induces the quotient mapping 
$\theta _M: \tilde M\vee \tilde M \to M\vee M$.
If both $M$ and $\tilde M$ 
satisfy the imposed above conditions on manifolds, then
this induces the embedding
$\bar \theta _M: Z^{\Upsilon ,a,c}(M,N)\hookrightarrow
Z^{\Upsilon ,a,c}(\tilde M,N)$. 
The identity mapping $id(x)=x$ for each $x\in \tilde M$
evidently satisfy Condition $(i)$.
If $\tilde f$ is the diffeomorphism of $\tilde M$
satisfying Condition $(i)$, then
applying $\tilde f^{-2}$ to both sides of the equality
we see, that it is satisfied for $\tilde f^{-1}$
with the same covering $\{ U_x: M_0 \} $.
If $\tilde f$ and $\tilde h$ are two diffeomorphisms of $\tilde M$,
then there exists $\{ U_x^h: x\in M_0^h \} $ such that
$\{ \tilde f^{-1} (\theta _M^{-1}(U_x^h)): x\in M_0^h \} $
is the locally finite covering of $\tilde M$. 
Two manifolds $M$ and $\tilde M$ are metrizabel, consequently, paracompact
(see Theorem 5.1.3 \cite{eng}).
Due to paracompactness of $M$ and $\tilde M$
there exists a locally finite covering
$\{ U_x^{h\circ f}: x\in M_0^{h\circ f} \} $ 
for which Condition $(i)$ is satisfied, since 
$\{ U_x^f\cap \theta _M\circ \tilde f^{-1} 
(\theta _M^{-1}(U_z^h)): x\in M_0^f, z\in M_0^h \} $
has a locally finite refinement.
This means, that $\bar \theta _M: Diff^{\infty }(M_{\bf R})\hookrightarrow
Diff^{\infty }(\tilde M_{\bf R})$ is the group embedding. 
In a complete uniform space $(X,{\sf U})$ for its subset $Z$ 
a uniform space $(Z,{\sf U}_Z)$ is complete if and only if
$Z$ is closed in $X$ relative to the topology induced by $\sf U$
(see Theorem 8.3.6 \cite{eng}).
Since both groups are complete and the uniformity of 
$Diff^{\infty }(\tilde M_{\bf R})$ induces the uniformity in 
$Diff^{\infty }(M_{\bf R})$ equivalent to its own, then
$\bar \theta _M( Diff^{\infty }(M_{\bf R}))$ is closed in
$Diff^{\infty }(\tilde M_{\bf R})$.
Considering the equivalence relation in
$Z^{\Upsilon ,a,c}(\tilde M,N)$ we get the loop monoids embedding
$\tilde \theta _M :(S^MN)_{\xi }\hookrightarrow
(S^{\tilde M}N)_{\xi }$ (see \S \S 2.1.4 and 2.1.5).
This produces with the help of A. Grothendieck construiction
the loop groups embedding (respecting their topological group structures)
$\tilde \theta _M:
(L^MN)_{\xi }\hookrightarrow (L^{\tilde M}N)_{\xi }.$
\par {\bf 2.2.1. Note.} For the diffeomorphism group we also consider
a compact complex manifold $M$.
For noncompact complex $M$, satisfying conditions of \S 2.1.1
and $(N1)$ the diffeomorphism group is considered
as consisting of diffeomorphisms $f$ of class $Y^{\Upsilon ,a,c}$
(see \S 2.1.5), that is, 
$(f_{i,j}-id_{i,j})\in Y^{\Upsilon ,a,c}(U_{i,j},\phi _i(U_i))$
for each $i, j$, $U_{i,j}$ is a domain of definition of 
$(f_{i,j}-id_{i,j})$ and then analogously to the real case
the diffeomorphism group $Diff^{\xi }(M)$
is defined, where $\xi =(\Upsilon ,a,c)$, 
$a=(a_1,a_2)$, $c=(c_1,c_2)$, 
$a_1 \le -1 $ and 
$c_1\le -1.$
This means that $Diff^{\xi }(M):=Y^{\xi ,id}(M,M)\cap Hom(M)$.
\par For investigations of stochastic processes on diffeomorphism groups
at first there are given below necessary definitions and statements 
on special kinds of diffeomorphism groups having Hilbert 
manifold structures.
\par{\bf 2.2.2. Remarks and definitions.} 
Let $M$ and $N$ be real manifolds on $\bf  R^n$ or $l_2$ 
and satisfying Conditions 2.2.(i-vi) \cite{lurim1}
or they may be also canonical closed submanifolds of that of
in \cite{lurim1}.
For a field ${\bf K}=\bf R$ or $\bf C$ let $l_{2,\delta }({\bf K})$
be a Hilbert space of vectors $x=(x^j: x^j\in {\bf K}, j\in {\bf N} )$
such that $\| x \| _{l_{2,\delta }}:=\{ \sum_{j=1}^{\infty }
|x^j|^2j^{2\delta } \} ^{1/2}<\infty $.
For $\delta =0$ we omit it as the index.
Let also $U$ be an open subset in $\bf R^m$ and $V$ be an open subset in
$\bf R^n$ or $l_2$ over $\bf R$, 
where $0\in U$ and $0\in V$ with $m$ and $n\in \bf N.$
By $H^{l, \theta }_{\beta ,\delta }(U,V)$
is denoted the following completion relative to the metric
$q^l_{\beta ,\delta }(f,g)$ of the family of all strongly
infinite differentiable functions 
$f,g: U \to V$ with $q^l_{\beta ,\delta }(f,\theta )<\infty $, where
$\theta \in C^{\infty }(U,V)$, $0\le l\in \bf Z$,
$\beta \in \bf R$, $\infty >\delta \ge 0$, $q^l_{\beta ,\delta }(f,g)
:=(\sum_{0 \le |\alpha |\le l} \| \bar m^{\alpha \delta }
<x> ^{\beta +|\alpha |}D_x^{\alpha }
(f(x)-g(x)) \|^2_{L^2})^{1/2}$, $L^2:=L^2(U,F)$
(for $F:=\bf R^n$ or $F=l_{2, \delta }=l_{2,\delta }({\bf R})$)
is the standard Hilbert space of all classes of equivalent 
measurable functions $h: U \to F$ for which there exists
$\| h\|_{L^2}:=(\int_U|h(x)|_F^2
\mu _m(dx)))^{1/2} <\infty $, $\mu _m$ denotes the Lebesgue measure
on $\bf R^m$. 
\par Let also $M$ and $N$ have finite atlases
such that $M$ be on
$X_M:=\bf R^m$ and $N$ on $X_N:=\bf R^n$ or $X_N:=l_2$, 
$\theta : M\hookrightarrow N$ be a $C^{\infty }$-mapping,
for example, embedding.
Then $H^{l,\theta }_{\beta ,\delta }(M,N)$ denotes the completion
of the family of all $C^{\infty }$-functions $g, f:
M \to N$ with $\kappa ^l_{\beta ,\delta }(f,\theta )< \infty $,
where the metric is given by the following formula
$\kappa ^l_{\beta ,\delta }(f,g)=(\sum_{i,j}[q^l_{
\beta ,\delta }(f_{i,j},
g_{i,j})]^2)^{1/2}$, where $f_{i,j}:=\psi _i\circ f\circ \phi _j^{-1}$ 
with domains $\phi _j(U_j)\cap \phi _j(f^{-1}(V_i))$,
$At(M):= \{ (U_i,\phi _i): i \} $ and $At(N):= \{ (V_j,\psi _j): j \} $
are atlases of $M$ and $N$, $U_i$ are open subsets in $M$ and 
$V_j$ are open subsets in $M$, $\phi _i: U_i\to X_M$ and 
$\psi _j: V_j\to X_N$ are homeomorphisms of $U_i$ on $\phi _i(U_i)$ and
$V_j$ on $\psi _j(V_j)$, respectively. 
Hilbert spaces $H^{l,\theta }_{\beta ,\delta }(U,F)$
and $H^l_{\beta ,0}(TM)$ are called weighted Sobolev spaces, 
where $H^l_{\beta ,\delta }(TM):=\{  f:M\to TM:$ $f \in H^l_{\beta ,\delta
}(M,TM),$ $\pi \circ
f(x)=x\mbox{ for each } x \in M \} $ with $\theta (x)=(x,0)\in T_xM$
for each $x\in M$. From the latter definition it follows, that
for such $f$ and $g$ there exists $\lim_{R\to \infty }q^l_{\beta ,\delta }
(f|_{U_R^c}, g|_{U_R^c})=0$, when $(U,\phi )$ is a chart Hilbertian 
at infinity, $U_R^c$ is an exterior of a ball of radius $R$ 
in $U$ with center in the fixed point $x_0$
relative to the distance function $d_M$ in $M$
induced by the Riemann metric $g$ (see \S 2.2(v) \cite{lurim1}). 
For $\beta =0$ or $\gamma =0$ we omit $\beta $
or $\gamma $ respectively in the notation 
$Dif^t_{\beta ,\gamma }(M):=H^{t,id}_{\beta ,\gamma }(M,M)\cap Hom(M)$ 
and $H^{l,\theta }_{\beta ,\gamma }$.
\par The uniform space $Dif^t_{\beta ,\gamma }(M)$ has the group structure
relative to the composition of diffeomorphisms and is called 
the diffeomorphism group, where $Hom(M)$ is the group of homeomorphisms 
of $M$. 
\par Each topologically adjoint space
$(H^l_{\beta }(TM))'=:H^{-l}_{-\beta }(TM)$ also is the Hilbert space
with the standard norm in
$H'$ such that $\| \zeta \|_{H'}=\sup_{\| f\| _H=1}| \zeta (f)|$.
\par  {\bf 2.3. Diffeomorphism groups of Gevrey-Sobolev classes of 
smoothness. Notes and definitions.} 
Let $U$ and $V$ be open subsets in the Euclidean space
$\bf R^k$ with $k\in \bf N$ or in the standard separable Hilbert 
space $l_2$ over $\bf R$,
$\theta : U\to V$ be a $C^{\infty }$-function (infinitely strongly 
differentiable), $\infty >\delta \ge 0$ be a parameter.
There exists the following metric space 
$H^{ \{ l \} ,\theta }_{ \{ \gamma \} ,\delta }(U,V)$
as the completion of a space of all functions
$Q:= \{ f:$ $f \in E^{\infty ,\theta }
_{\infty ,\delta }(U,V),$ $\mbox{ there exists }n \in {\bf N}
\mbox{ such that }$
$supp (f) \subset U\cap {\bf R^n}, d_{ \{ l\}, \{ \gamma \},\delta }
(f,\theta )<\infty \} $
relative to the given below metric
$d_{ \{ l\}, \{ \gamma \},\delta }:$ 
$$(i)\mbox{ }d_{ \{ l\}, \{ \gamma \},\delta }(f,g):=
\sup_{x \in U}(\sum_{n=1}^{\infty }(\bar \rho ^l_{\gamma ,n,\delta }
(f,g))^2)^{1/2}<\infty $$ and 
\par $\lim_{R \to \infty }d_{\{ l\},\{ \gamma \},\delta }(f|_{U_R^c},
g|_{U_R^c})=0$, when $U$ is a chart Euclidean or Hilbertian 
correspondingly at infinity,
$f$ as an argument in $\bar \rho ^l_{\gamma ,n,\delta }$ 
is taken with the restriction on
$U\cap \bf R^n$, that is,  $f|_{U\cap \bf R^n}: U\cap {\bf R^n}\to
f(U)\subset V$ (see also \S \S 2.1-2.5 \cite{lurim1} and \cite{lurim2}
about $E^{t,\theta }_{\beta ,\gamma }$), 
$\bar \rho ^l_{\gamma ,n,\delta }(f,id)^2$ 
$:=\omega _n^2(\kappa ^l_{\gamma ,\delta }(f|_{(U\cap {\bf R^n})},
id|_{(U \cap {\bf R^n})})^2 -$
$\kappa ^{l(n-1)}_{\gamma (n-1),\delta }(f|_{(U\cap {\bf R^{n-1}})},
id|_{(U \cap {\bf R^{n-1}})})^2)$ 
for each $n>1$ and \\
$\bar \rho ^l_{\gamma ,1,\delta }(f,id)
:=\omega _1(\kappa ^l_{\gamma ,\delta }(f|_{(U\cap {\bf R^1})},
id|_{(U \cap {\bf R^1})})$ 
with $q^l_{\gamma ,\delta }$
and the corresponding terms $\kappa ^l_{\gamma ,\delta }$ 
from \S 2.2.2,
$l=l(n)>n+5$, $\gamma =\gamma (n)$ and $l(n+1)\ge l(n)$
for each $n$, $l(n)\ge [t]+sign \{ t \}+[n/2]+3$, $\gamma
(n)\ge \beta +sign \{ t \} +[n/2]+7/2$, $\omega _{n+1}\ge n\omega _n\ge 1$. 
Moreover, $\bar \rho ^l_{\gamma ,n,\delta }
(f,id)(x^{n+1},x^{n+2},...)\ge 0$ is the metric by variables
$x^1,...,x^n$ in $H^l_{\gamma ,\delta }(U\cap
{\bf R^n}, V)$ for $f$ as a function by 
$(x^1,...,x^n)$ such that $\bar \rho ^l_{\gamma ,n,\delta }$
depends on parameters $(x^j: j>n)$. The index $\theta $ 
is omitted when $\theta =0$.
The series in $(i)$ terminates $n\le k$, when $k\in \bf N$.
\par Let for $M$ connecting mappings of charts be such that
$(\phi _j\circ \phi _i^{-1}-id_{i,j})
\in H^{ \{ l' \} }_{ \{ \gamma ' \} ,\chi }(U_{i,j},l_2)$ for each $U_i
\cap U_j\ne \emptyset $ and the Riemann metric $g$ be of class
of smoothness $H^{
\{ l' \} }_{ \{ \gamma ' \} ,\chi }$, where subsets $U_{i,j}$ are open in
$\bf R^k$ or in $l_2$ correspondingly
domains of $\phi _j \circ \phi _i^{-1}$, $l'(n)\ge l(n)+2$,
$\gamma '(n)\ge \gamma (n)$ for each $n$, $\infty >\chi \ge \delta $, 
submanifolds $\{ M_k: k=k(n), n \in {\bf N} \}$
are the same as in Lemma 3.2 \cite{lurim1}. 
Let $N$ be some manifold satisfying analogous conditions as $M$.
Then there exists the following uniform space
$H^{ \{ l \} ,\theta }_{ \{ \gamma \} ,\delta ,\eta }(M,N):=
\{ f\in E^{\infty ,\theta }_{\infty ,\delta }(M,N)|
(f_{i,j}-\theta _{i,j})\in H^{ \{ l \} ,\theta }_{ \{ \gamma \} ,\delta }
(U_{i,j},l_2)$ $\mbox{for each charts }$
$\{ U_i,\phi _i \}$ and $\{ U_j,\phi _j \}$ $\mbox{ with }U_i\cap U_j
\ne \emptyset , \chi _{\{ l\},\{ \gamma \},\delta ,\eta }
(f,\theta )<\infty $ $\mbox{and }\lim_{R\to \infty }
\chi _{\{ l\},\{ \gamma \},\delta ,\eta }(f|_{M^c_R},\theta 
|_{M^c_R})=0 \} $ and there exists the corresponding diffeomorphism group
$Di^{ \{ l \} }_{ \{ \gamma \} , \delta ,\eta }(M):=$ $ \{ f : f\in Hom(M),$
$f^{-1}$ $\mbox{ and }$ $f\in
H^{ \{ l \} ,id}_{ \{ \gamma \} ,\delta ,\eta }(M,M) \} $ 
with its topology given by the following left-invariant metric
$\chi _{\{ l\}, \{ \gamma \},\delta ,\eta }(f,g):=\chi _{ \{ l\},
\{ \gamma \}, \delta ,\eta }(g^{-1}f,id)$,
$$(ii)\mbox{ }\chi _{\{ l\}, \{ \gamma \}, \delta ,\eta }
(f,g):=(\sum_{i,j}(d_{ \{ l\}, \{ \gamma \}, \delta }(f_{i,j},
g_{i,j})i^{\eta }j^{\eta })^2)^{1/2}<\infty ,$$  
$g_{i,j}(x)\in l_2$ and $f_{i,j}(x) \in l_2$,
$\phi _i(U_i)\subset l_2$, $U_{i,j}=U_{i,j}(x^{n+1},x^{n+2},...)
\subset l_2$ are domains of $f_{i,j}$ by variables 
$x^1,...,x^n$ for chosen $(x^j: j>n)$ due to existing foliations in $M$, 
$U_{i,j} \subset {\bf R^n}\hookrightarrow l_2$, when $(x^j: j>n)$
are fixed and $U_{i,j}$ is a domain in $\bf R^n$ by variables
$(x^1,...,x^n)$, where $\infty >\eta \ge 0$.
\par In particular, for the finite dimensional manifold
$M_n$ the group $Di^{ \{ l\} }_{ \{ \gamma \} ,\delta ,\eta }(M_n)$ 
is isomorphic to the diffeomorphism group 
$Dif^l_{\gamma ,\delta }(M_n)$
of the weighted Sobolev class of smoothness
$H^l_{\gamma ,\delta }$ with $l=l(n)$, $\gamma =\gamma (n)$,
where $n=dim_{\bf R}(M_n)<\infty $.
\par {\bf 2.4. Remarks.} Let two sequences be given
$\{ l \}:=\{ l(n): n \in {\bf N} \}\subset
\bf Z$ and $\{ \gamma \} :=\{ \gamma (n): n \in {\bf N} \}\subset
\bf R$, where $M$ and $\{ M_k: k=k(n), n \in {\bf N} \} $ are
the same as in \S \S 2.2 and 2.3. Then there exists the following space 
$H^{ \{ l \} ,\theta }_{ \{ \gamma \} ,\delta ,\eta }(M,TN)$.
By $H^{ \{ l \} ,\theta }_{ \{ \gamma \} ,\delta ,\eta }(M|TN)$
it is denoted its subspace of functions $f: M\to TN$
with $\pi _N(f(x))=\theta (x)$ for each $x\in M$,
where $\pi _N: TN\to N$ is the natural projection,
that is, each such 
$f$ is a vector field along $\theta $, $\theta : M\to N$ is a
fixed $C^{\infty }$-mapping. For $M=N$
and $\theta =id$ the metric space
$H^{ \{ l \} ,\theta }_{ \{ \gamma \} ,\delta ,\eta }(M|TM)$
is denoted by
$H^{ \{ l \} }_{ \{ \gamma \} ,\delta ,\eta }(TM)$. Spaces
$H^{ \{ l \} ,id}_{ \{ \gamma \} ,\delta ,\eta }(M_k|TN)$
and $H^{ \{ l \} }_{ \{ \gamma \} ,\delta ,\eta }(TM)$
are Banach spaces with the norms 
$\| f\|_{ \{ l\}, \{ \gamma \} ,\delta ,\eta }
:=\chi _{\{ l\}, \{ \gamma \}, \delta ,\eta }(f,f_0)$
denoted by the same symbol, where
$f_0(x)=(x,0)$ and $pr_2f_0(x)=0$
for each $x\in M$. This definition can be spread on the case
$l=l(n)<0$, if take $\sup
_{\| \tau \| =1}| <x>_m^{| \alpha |- \gamma (m)}$ 
$(D_x^{\alpha }\tau _{i,j},
[\zeta _{i,j}-\xi _{i,j}])_{L^2(U_{i,j,m},l_{2,\delta })}|$
instead of $\| <x>_m^{\gamma (m)+|
\alpha |}$ $D_x^{\alpha }(\zeta _{i,j}-\xi _{i,j})(x) \| _{L^2(
U_{i,j,m}, l_{2,\delta })}$, where
$\tau \in H_{-\gamma } ^{-l}(M_k|TN)$, $<x>_m=(1+\sum_{i=1}^m(x^i)^2)^{1/2}$,
$U_{i,j,m}=U_{i,j,m}(x^{m+1},x^{m+2},...)$ denotes the 
domain of the function $\zeta _{i,j}$
by $x^1,...,x^m$ for chosen $(x^j:$ $j>m)$, 
$\| \zeta \| _k(x)$ are functions by variables $(x^i: i>k)$.
Further the traditional notation is used:  
$sign( \epsilon )=1$ for $\epsilon >0$, $sign(\epsilon )=-1$
for $\epsilon <0$, $sign (0)=0$, $\{ t \}=t-[t]\ge 0$.
\par  {\bf 2.5. Lemma.} {\it Let the manifold $M$
and the spaces $E^t_{\beta ,\delta }(TM)$
and $H^{ \{ l \} }_{\{ \gamma \} ,\delta ,\eta }(TM)$ 
be the same as in \S \S 2.2-2.4 with $l(k)\ge
[t]+[k/2]+3+sign \{ t \}$, $\gamma (k)\ge \beta +[k/2]+7/2+sign\{ t \}$.
Then there exist constants $C>0$ and $C_n>1$ for each $n$
such that $\| \zeta \| _{E^t_{\beta ,\delta }(TM)} \le 
C \| \zeta \| _{ \{ l \}, \{ \gamma \} ,\delta ,0}$ for each
$\zeta \in H^{ \{ l \} }_{ \{ \gamma \} ,\delta ,0}(TM)$, 
moreover, there can be chosen $\omega _n\ge C_n$, $C_{n+1}\ge k(n+1)
(k(n+1)-1)...(k(n)+1)C_n$
for each $n$ such that the following inequality be valid:
$\| \xi \| _{C^{l'(k)}_{\gamma '(k)}(TM_k)}$
$ \le C_n \| \xi \| _{H^{l(k)}
_{\gamma (k)}(TM_k)}$ for each $k=k(n)$, $l'(k)=l(k)-[k/2]-1$,
$\gamma '(k)=
\gamma (k)-[k/2]-1$ for each $\xi \in H^{l(k)}_{\gamma (k)}(TM_k)$.}
\par  {\bf Proof.} In view of theorems from ~\cite{tri} 
and the inequality $\int_{\bf R^m}
<x>_m^{-m-1}dx< \infty $ (for $<x>_m$ taken in $\bf R^m$ with
$x \in \bf R^m$) there exists the embedding $H^{l(n)}_{\gamma (n)}(TM_n)$
$\hookrightarrow C^{l'(n)}_{\gamma '(n)}(TM_n)$ for each $n$, 
since $2([n/2]+1)
\ge n+1$.  Moreover, due to results of \S III.6 \cite{miha}
there exists a constant $C_n>0$ for each $k=k(n)$, 
$n \in \bf N$ such that $\| \xi \| _{C^{l'(k)}_{\gamma '(k)}(TM_k)}$
$ \le C_n \| \xi \| _{H^{l(k)}
_{\gamma (k)}(TM_k)}$ for each $\xi \in H^{l(k)}_{\gamma (k)}(TM_k)$.
Then $D^{\alpha }f(x^1,...,x^n,...)-D^{\alpha }f(y^1,...,y^n,...)=$
$$\sum_{n=0}^{\infty }(D^{\alpha }f(y^1,...,y^{n-1},x^n,...)-D^{\alpha }
f(y^1,...,y^n,x^{n+1},...))$$ 
for each $f \in H^{ \{ l \} }_{ \{ \gamma \} ,\delta }(TM)$
in local coordinates, where $f(y^1,...,y^{n-1},x^n,...)=f(x^1,x^2,...,
x^n,...)$, if $n=0$; $\alpha =(\alpha ^1,...,\alpha ^m)$, $m \in \bf N$,
$\alpha ^i \in {\bf N_o}:=\{ 0,1,2,... \} $. 
Hence for $x^n< y^n$ the following inequlity is satisfied: \\
$|D^{\alpha }f(y^1,...,y^{n-1},
x^n,x^{n+1},...)-D^{\alpha }f(y^1,...,y^n,x^{n+1},...)| _{l_{2,\delta }}
\bar m^{\alpha \delta }\le $
$$[\int_{\phi _j(U_j\cap M_k)\ni z:=
(y^1,...,y^{n-1},z^n), x^n\le z^n\le y^n}
|D^{\alpha }\partial f(y^1,...,y^{n-1},z^n,x^{n+1},...)/\partial z^n
|_{l_{2,\delta }}dz^n]\bar m^{\alpha \delta }\le $$
$$C^1\int \int_{\phi _j(U_j\cap M_{k(n+1)}) 
\ni z:=(y^1,...,y^{n-1},z^n, z^{n+1}), x^n\le z^n\le y^n}
sup_{x \in M}(\| f \| _{H^{l(k(n+8))}_{ \gamma (k(n+8))}(M_{k(n+8)}|TM)}$$
$<z>_{n+1}^{-5/2})dz^ndz^{n+1}(n+1)^{-2}
\le C' \| f \| _{ \{ l \}, \{ \gamma \} ,\delta ,0} \times (n+1)^{-2}$, \\
when $| \alpha |=\alpha ^1+...+\alpha ^m \le l(k)$,
$k=k(n)\ge n$, $m \le n$,
where $C^1=const >0$ and $C'=const >0$ are constants not depending on
$n$ and $k$; $x,$ $y$ and $(y^1,..,y^n,x^{n+1},x^{n+2},...)\in \phi _j(U_j)$
for each $n\in \bf N$. 
This is possible due to local convexity of the subset 
$\phi _j(U_j)\subset l_2$.
Therefore, $H^{ \{ l \} }_{ \{
\gamma \} ,\delta ,0}(TM) \subset E^t_{\beta ,\delta }(TM)$ and $\| f \|
_{E^t_{\beta ,\delta }(TM)}\le C \| f \| _{ \{ l \}, \{ \gamma \} ,
\delta ,0}$ for each
$f \in H^{ \{ l \} }_{ \{ \gamma \} ,\delta ,0}(TM)$, moreover, $C=C'
\sum_{n=1}^{\infty }
n^{-2} <\infty $, since $\sup _{x \in M}\sum_{j=1}^{\infty }g_j(x)
\le \sum_{j=1}^{\infty }\sup_{x \in M}g_j(x)$ for each function 
$g: M \to [0, \infty )$ and
$\lim_{R\to \infty }\| f|_{M_R^c} \| _{E^t_{\beta ,\delta }(TM)}$
$\le C\times \lim_{R\to
\infty }\| $ $f|_{M_R^c}\|_{\{ l \}, \{ \gamma \} ,\delta ,0}=0$.
\par The space 
$E^t_{\beta ,\delta }(TM)\cap H^{ \{ l \} }_{ \{ \gamma \} ,\delta,0}
(TM)$ contains the corresponding cylindrical functions
$\zeta $, in particular with
$supp(\zeta ) \subset U_j\cap M_n$ for some $j \in \bf N$ and
$k=k(n)$, $n \in \bf N$. The linear span of the family $\sf K$
over the field $\bf R$ of such functions $\zeta $ is dense in 
$E^t_{\beta ,\delta }(TM)$ and in $H^{ \{ l \} }_{ \{ \gamma \} ,\delta ,0}
(TM)$ due to the Stone-Weierstrass theorem,
consequently, $H^{ \{ l \} }_{ \{ \gamma \},\delta ,0}(TM)$ 
is dense in $E^t_{\beta ,\delta }(TM)$,
since $\partial f/\partial x^{n+1}=0$
for cylindrical functions $f$ independent from $x^{n+j}$ for $j>0$.
\par {\bf 2.6.1. Note.}
For the diffeomorphism group $Diff^t_{\beta ,\gamma }(\tilde M)$
of a Banach manifold $\tilde M$ let $M$ be a dense Hilbert submanifold
in $\tilde M$ as in \cite{lurim1,lurim2}.
\par  {\bf 2.6.2. Lemma.} {\it Let 
$Di^{ \{ l \} }_{ \{ \gamma \} , \delta ,\eta }(M)$ and $M$ 
be the same as in \S 2.3 with values of parameters $C_n$ from Lemma 2.5 
for given $l(k)$, $\gamma (k)$ and $k=k(n)$ with $\omega _n=l(k(n))!C_n$, 
then $Di^{ \{ l \} }_{ \{ \gamma \} , \delta ,\eta }(M)$ is the 
separable metrizable topological 
group dense in $Diff^t_{\beta ,\delta }(\tilde M)$.
\par In the case of the complex manifold $M$ 
the group $Diff^{\xi }(M)$ (see \S \S 2.1.5 and 2.2.1) is the separable
metrizable topological group.}
\par  {\bf Proof.} Consider at first the real case.
From the results of the paper \cite{omo}
it follows that the uniform space
$Di^{ \{ l \} }_{ \{ \gamma \} , \delta ,\eta }(M_k)$ 
is the topological group for each finite dimensional submanifold
$M_k$, since $l(k)>k+5$ and $dim_{\bf R}M_k=k$.
The minimal algebraic group
$G_0:=gr(Q)$ generated by the family $Q:=\{ f: f\in E^{ \{ l\} ,id}_{
\{ \gamma \} ,\delta }(U,V)$ $\mbox{for all possible pairs of charts }
U_i$ $\mbox{and}$ $U_j$ $\mbox{with}$ $U=\phi _i(U_i)$ $\mbox{and}$
$V=\phi _j(U_j),$ $supp(f)\subset U\cap \bf R^n,$ 
$f\in Hom(M),$ $dim_{\bf R}M\ge n \in {\bf N} \}$
is dense in $Di^{ \{ l \} }_{ \{ \gamma \} , \delta ,\eta }(M)$ 
and in $Diff^t_{\beta ,\delta }(M)$ due to the Stone-Weierstrass theorem, 
since the union $\bigcup_kM_k$ is dense in $M$, where
$supp(f):=cl\{ x\in M: f(x)\ne x \}$, $cl(B)$ denotes the closure
of a subset $B$ in $M$. 
Therefore, $Di^{ \{ l \} }_{ \{ \gamma \} , \delta ,\eta }(M)$
and $Diff^t_{\beta ,\delta }(\tilde M)$ are separable.
It remains to verify that
$Di^{ \{ l \} }_{ \{ \gamma \} , \delta ,\eta }(M)$
is the topological group. For it we shall use 
Lemma 2.5. For $a>0$ and $k\ge 1$ using integration by parts 
formula we get the following equality $\int^{\infty }_{-\infty }(
a^2+x^2)^{-(k+2)/2}dx=((k-1)/(ka^2))\int^{\infty }_{-\infty }
(a^2+x^2)^{k/2}dx$, which takes into account the weight multipliers.
Let $f, g\in Di^{ \{ l\} }_{ \{ \gamma \}, \delta ,\eta }(V)$
for an open subset $V=\phi _j(U_j)\subset l_2$ and
$\chi _{ \{ l\}, \{ \gamma \},\delta ,\eta }(f,id)<1/2$
and $\chi _{ \{ l\}, \{ \gamma \},\delta ,\eta }(g,id)<\infty $, then
$\bar \rho ^l_{\gamma ,n,\delta }(g^{-1}\circ f,
id) $ $\le C_{l,n,\gamma ,\delta }(\bar \rho ^{4l}_{\gamma ,n,\delta }
(f,id)+\bar \rho ^{4l}_{\gamma ,n,\delta }(g,id))$,
where $0<C_{l,n,\gamma ,\delta }\le 1$ is a constant dependent on 
$l, n, \gamma $ and independent from $f$ and $g$.
For the Bell polynomials there is the following inequality
$Y_n(1,...,1) \le n!e^n$ for each $n$ and $Y_n(F/2,...,F/(n+1))
\le (2n)!e^n$ for $F^p:=F_p=(n+p)_p:=(n+p)...(n+2)(n+1)$ (see
Chapter 5 in ~\cite{rio} and Theorem 2.5 in ~\cite{ave}). 
The Bell polynomials are given by the following formula
$Y_n(fg_1,...,fg_n):=\sum_{\pi (n)}(n!f_k/(k_1!...k_n!))(g_1/1!)^{k_1}
...(g_n/n!)^{k_n}$, where the sum is by all partitions
$\pi (n)$ of the number $n$, this partition is denoted by $1^{k_1}2^{k_2}...
n^{k_n}$ such that $k_1+2k_2+...+nk_n=n$ and $k_i$ is a number of terms
equal to $i$, the total number of terms in the partition is equal to
$k=k(\pi )=k_1+....+k_n$, $f^k:=f_k$ in the Blissar calculus notation.
For each $n\in \bf N$, $l=l(n)$ and $\gamma =\gamma (n)$
the following inequality is satisfied:
$\bar \rho ^l_{\gamma ,n,\delta }(f\circ g,id)
\le Y_l(\bar f\bar g_1,...,\bar f\bar g_m)$, 
$\bar \rho ^{l+1}_{\gamma ,n,\delta }(f^{-1},id)\le (3/2)Y_l(Fp_1/2,
...,Fp_m/(m+1))$, where $\bar f^m:=\bar f_m=\bar \rho ^m_{\gamma ,n,
\delta }(f,id)$ and $F^k:=F_k=(n+k)_k$, $(n)_j:=n(n-1)...(n-j+1)$,
$p_k=-\bar f_{k+1}(3/2)^{k+1}$.
Then $\sum_{n=1}^{\infty }(2l(k(n)))!e^{l(k(n))}b^{4l(k(n))}
(l(k(n))!)[(4l(k(n)))!]^{-1}<\infty $ for each $0<b<\infty $.
Hence due the Cauchy-Schwarz-Bunyakovskii inequality 
and the condition $C_n^2>C_n$
for each $n$ we get: $f\circ g$ and $f^{-1}\in Di^{ \{ l\} }_{ \{ \gamma \},
\delta ,\eta }(M)$ for each $f$ and $g\in 
Di^{ \{ l\} }_{ \{ \gamma \},\delta ,\eta }(M)$, moreover, the operations 
of composition an inversion are continuous.
\par The base of neighborhoods of $id$ in
$Di^{ \{ l \} }_{ \{ \gamma \} , \delta ,\eta }(M)$
is countable, hence this group is metrizable, moreover, 
a metric can be chosen left-invariant due to Theorem 8.3 \cite{hew}.
The case $Diff^{\xi }(M)$ for the complex manifold $M$ is analogous.
\par {\bf 2.7. Lemma.} {\it Let 
$G':=Di^{ \{ l" \} }_{ \{ \gamma " \} ,\delta ",\eta "}(M)$ 
be a subgroup of
$G:=Di^{ \{ l \} }_{ \{ \gamma \} ,\delta ,\eta }(M)$ 
such that $m(n)>n/2$,
$l"(n)=l(n)+m(n)n$, $\gamma "(n)=\max (\gamma (n)-m(n)n,0)$
for each $n$, $inf-\lim_{n\to \infty }m(n)/n=c>1/2$, 
$\delta ">\delta +1/2$, $\infty >\eta ">\eta +1/2,$ $\eta \ge 0$ 
(see \S 2.3). 
Let also $G':=Diff^{\xi '}(M)$ be a subgroup of
$G=Diff^{\xi }(M)$ with either ${a'}_1<a_1$ 
and ${c'}_1<c_1$ or ${a'}_1=a_1$ and ${a'}_2<a_2-1$ 
and ${c'}_1=c_1$ and ${c'}_2<c_2-1$ for the complex manifold $M$ 
(see \S \S 2.1.5 and 2.2.1).
Then there exists a Hilbert-Schmidt operator of
embedding $J : Y'\hookrightarrow Y$, where
$Y:=T_eG$ and $Y':=T_eG'$ are tangent Hilbert spaces.}
\par  {\bf Proof.} Consider at first the real case.
The natural embeddings 
$\theta _k$ of the Hilbert spaces $H^{l(k)-m(k)k,b(k)}_{
\gamma (k)+m(k)k,\delta }(M_k,{\bf R})$
into $H^{l(k),b(k)}_{\gamma (k)k,\delta }(M_k,l_{2,\delta +1+\epsilon })$
are Hilbert-Schmidt operators for each $k=k(n)$, $n\in \bf N$. 
For each chart $(U_j,\phi _j)$ there are 
linearly independent functions
$x^me_l<x>^{\zeta }_n/m!=:f_{m,l,n}(x)$, where $\{ e_l: l \in {\bf N} \}
\subset l_2$ is the standard orthonormal basis in $l_2$, 
$x^m:=x_1^{m_1}...x_n^{m_n}$,
$m!=m_1!...m_n!$, $<x>_n=$ $(1+\sum_{i=1}^n(x^i)^2)^{1/2}$, $n \in \bf N$,
$\zeta (n)=\zeta \in \bf R$.
The linear span over $\bf R$ of the family of all such functions
$f(x)$ is dense in $Y$.
Moreover, $D^{\alpha }f(x)=e_l\sum {\alpha \choose \beta }(D^{\beta }x^m/m!)$
$(D^{\alpha -\beta }<x>^{\zeta }_n)$, where 
$D^{\alpha }=\partial _1^{\alpha ^1}...\partial _n^{\alpha ^n},$
$\partial _i=\partial /\partial x_i$,
$\alpha =(\alpha ^1,...,
\alpha ^n)$, ${\alpha \choose \beta }$
$=\prod _{i=1}^n{{\alpha ^i} \choose {\beta ^i}}$,
$0\le \alpha ^i\in \bf Z$, $\lim_{n \to \infty }
q^n/n!=0$ for each $\infty >q>0$, $\sum_{j,l,n=1}^{\infty }
\sum_{|m|
\ge m(n), m}[jln^nm_1...m_n]^{-(1+2\epsilon )}<\infty $
for each $0<\epsilon <\min (c-1/2, \eta "- \eta -1/2, \delta "-\delta -1/2)$, 
where $m=(m_1,...,m_n)$, $|m|:=m_1+...+m_n$, $0\le m_i \in \bf Z$.
Hence due to \S \S 2.3 and 2.4 the embedding $J$ is the Hilbert-Schmidt 
operator.
\par In the complex case we use the convergence of the series \\
$\sum_{j=1}^{\infty }\sum_{n=1}^{\infty}(j!)^{{a'}_1-a_1}
(n!)^{{c'}_1-c_1}<\infty $ and 
$\sum_{j=1}^{\infty }\sum_{n=1}^{\infty}j^{{a'}_2-a_2}
n^{{c'}_2-c_2}<\infty $.
\par For the construction of Wiener processes on loop and diffeomorphism 
groups the existence of uniform atlases for them as manifolds
is necessary, that is given by the following proposition.
\par  {\bf 2.8. Theorem.}{\it Let the diffeomorphism group
$G:=Di^{ \{ l \} }_{ \{ \gamma \} , \delta ,\eta }(M)$
be the same as in \S 2.3 or $G:=Diff^{\xi }(M)$ as in  \S 2.2.1. Then
$$(i) \mbox{for each }H^{ \{l \} ,id }_{ \{ \gamma \} ,\delta ,\eta }
(M,TM)\mbox{-vector field } V \mbox{ its
flow }\eta _t$$ 
$\mbox{ is a one-parameter subgroup of }
Di^{ \{ l \} }_{ \{ \gamma \} , \delta ,\eta }(M)$,
$\mbox{ the curve }t\mapsto \eta _t\mbox{ is of class }C^1,$ 
$\mbox{ the mapping}$
$\tilde Exp: \tilde T_e Di^{ \{ l \} }_{ \{ \gamma \} , \delta ,\eta }(M)\to
Di^{ \{ l \} }_{ \{ \gamma \} , \delta ,\eta }(M)$,
$\mbox{ is continuous }$ and defined on the neighbourhood 
$\tilde T_e Di^{ \{ l \} }_{ \{ \gamma \} , \delta ,\eta }(M)$
of
the zero section in
$T_eDi^{ \{ l \} }_{ \{ \gamma \} , \delta ,\eta }(M)$,
$V\mapsto \eta _1$;
$$(ii)\mbox{ }T_fDi^{ \{ l \} }_{ \{ \gamma \} , \delta ,\eta }(M),
=\{ V \in H^{ \{l \} ,id }_{ \{ \gamma \} ,\delta ,\eta }
(M,TM)| \pi \circ V=f \};$$
$$(iii)\mbox{ }(V,W)=\int_Mg_{f(x)}(V_x,W_x)\mu (dx)$$
is a weak Riemannian structure on a Hilbert 
manifold $Di^{ \{ l \} }_{ \{ \gamma \} , \delta ,\eta }(M)$,
where
$\mu $ is a measure induced on $M$ by $\phi _j$ and a Gaussian measure
with zero mean value on $l_2$ produced by an injective self-adjoint
operator $Q: l_2 \to l_2$ of trace class, $0< \mu (M) <\infty $;
$$(iv)\mbox{ the Levi-Civita connection }\nabla \mbox{ on } M
\mbox{ induces the Levi-
Civita connection }$$ $\hat \nabla $ on 
$Di^{ \{ l \} }_{ \{ \gamma \} , \delta ,\eta }(M)$;
$$(v)\mbox{ } \tilde E: T Di^{ \{ l \} }_{ \{ \gamma \} , \delta ,\eta }(M)
\to Di^{ \{ l \} }_{ \{ \gamma \} , \delta ,\eta }(M)
\mbox{ is
defined by }$$ $\tilde E_{\eta }(V)=exp_{\eta (x)} \circ V_{\eta }$
on a neighbourhood $\bar V$ of the zero section
in $T_{\eta }Di^{ \{ l \} }_{ \{ \gamma \} , \delta ,\eta }(M)$
and is a $H^{ \{l \} ,id }_{ \{ \gamma \} ,\delta ,\eta }$-mapping 
by $V$ onto a
neighbourhood $W_{\eta }=W_{id}\circ \eta $ of $\eta
\in Di^{ \{ l \} }_{ \{ \gamma \} , \delta ,\eta }(M)$; 
$\tilde E$ is the uniform isomorphism
of uniform spaces $\bar V$ and $W$.
Analogous statements are true for $Diff^{\xi }(M)$
with the class of smoothness $Y^{\xi ,id}$
instead of $H^{ \{l \} ,id }_{ \{ \gamma \} ,\delta ,\eta }$.}
\par   {\bf Proof.} Consider at first the real case.
Then we have that $T_f
H^{ \{l \} ,\theta }_{ \{ \gamma \} ,\delta ,\eta }(M,N)
= \{ g \in H^{ \{l \} ,\theta }_{ \{ \gamma \} ,\delta ,\eta }(M,TN):$
$ \pi _N\circ g=f \} $,
where $\pi _N: TN\to N$ is the canonical projection.
Therefore, $TH^{ \{l \} ,\theta }_{ \{ \gamma \} ,\delta ,\eta }(M,N)
=H^{ \{l \} ,\theta }_{ \{ \gamma \} ,\delta ,\eta }(M,TN)
=\bigcup_fT_fH^{ \{l \} ,\theta }_{ \{ \gamma \} ,\delta ,\eta }(M,N)$ 
and the following mapping $w_{exp}:T_f
H^{ \{l \} ,\theta }_{ \{ \gamma \} ,\delta ,\eta }(M,N) \to
H^{ \{l \} ,\theta }_{ \{ \gamma \} ,\delta ,\eta }(M,N)$, 
$w_{exp}(g)=exp \circ g$ gives
charts for $H^{ \{l \} ,\theta }_{ \{ \gamma \} ,\delta ,\eta }(M,N)$, 
since  $TN$ has an atlas of class 
$H^{ \{ l'(n)-1: n \} }_{ \{ \gamma ' (n)+1 : n \} , \chi }$
In view of Theorem 5 about differential equations
on Banach manifolds in \S 4.2 \cite{lan} a
vector field $V$ of class $H^{ \{l \} ,\theta }_{ \{ \gamma \} ,
\delta ,\eta }$ on $M$ defines a flow $\eta _t$
of such class, that is $d \eta _t/dt=V \circ \eta _t$ and
$\eta _0=e$. From the proofs of Theorem 3.1 and Lemmas 3.2,
3.3 in \cite{ebi} we get that $\eta _t$ is a one-parameter subgroup of
$Di^{ \{ l \} }_{ \{ \gamma \} , \delta ,\eta }(M)$, 
the curve $t \mapsto \eta _t$ is of class $C^1$, 
the map $\tilde Exp: T_eDi^{ \{ l \} }_{ \{ \gamma \} , \delta ,\eta }(M)
\to Di^{ \{ l \} }_{ \{ \gamma \} , \delta ,\eta }(M)$ defined by $V \mapsto
\eta _1$ is continuous.
\par The curves of the form $t \mapsto \tilde E(tV)$ are geodesics for $V
\in T_{\eta }Di^{ \{ l \} }_{ \{ \gamma \} , \delta ,\eta }(M)$ 
such that $d \tilde E(tV)/dt$ is the map
$m \mapsto d(exp(tV(m))/dt=\gamma '_m(t)$ for each $m\in M$, 
where $\gamma _m(t)$ is the
geodesic on $M$, $\gamma _m(0)=\eta (m)$, $\gamma _m'(0)=V(m)$.
Indeed, this follows from the existence of solutions of corresponding
differential equations in the Hilbert space 
$H^{ \{l \} ,\eta }_{ \{ \gamma \} ,\delta ,\eta }(M|TM)$,
then as in the proof of Theorem 9.1 \cite{ebi}.
\par From the definition of $\mu $ it follows that for each $x \in M$
there exists an open neighbourhood $Y \ni x$ such that $\mu (Y)>0$
\cite{sko}. 
Since $t\ge 1$, the scalar product $(iii)$ gives a weaker topology
than the initial $H^{ \{l \} }_{ \{ \gamma \} ,\delta ,\eta }$.
\par Then the right multiplication $\alpha _h(f)=f \circ h$, 
$f \to f \circ h$
is of class $C^{\infty }$ on 
$Di^{ \{ l \} }_{ \{ \gamma \} , \delta ,\eta }(M)$ 
for each $h \in Di^{ \{ l \} }_{ \{ \gamma \} , \delta ,\eta }(M)$. 
Moreover, $Di^{ \{ l \} }_{ \{ \gamma \} , \delta ,\eta }(M)$
acts on itself freely from the right, hence we have
the following principal vector bundle $\tilde \pi :
T Di^{ \{ l \} }_{ \{ \gamma \} , \delta ,\eta }(M) \to 
Di^{ \{ l \} }_{ \{ \gamma \} , \delta ,\eta }(M)$
with the canonical projection $\tilde \pi $.
\par Analogously to \cite{ebi,lurim1} we get the connection $\hat \nabla
=\nabla \circ h$ on 
$Di^{ \{ l \} }_{ \{ \gamma \} , \delta ,\eta }(M)$. 
If $\nabla $ is torsion-free then $\hat \nabla $
is also torsion-free. From this it follows that the existence of $\tilde E$
and $Di^{ \{ l \} }_{ \{ \gamma \} , \delta ,\eta }(M)$
is the Hilbert manifold of class
$H^{ \{ l'(n)-1: n \} }_{ \{ \gamma '(n)+1: n \} ,\chi ,\eta }$, 
since $exp$ for $M$ is of class
$H^{ \{ l'(n)-1: n \} }_{ \{ \gamma '(n)+1: n \} ,\delta }$,
$f\to f\circ h$ is a
$C^{\infty }$ map with the derivative $\alpha _h: 
H^{ \{l \} ,\eta }_{ \{ \gamma \} ,\delta ,\eta }(M',TN)
\to H^{ \{l \} ,\eta }_{ \{ \gamma \} ,\delta ,\eta }(M,TN)$ 
whilst $h \in H^{ \{l \} ,\eta }_{ \{ \gamma \} ,\delta ,\eta }(M,M'),$ 
\par $(vi)$ $\tilde E_h(\hat V):=
exp_{h(x)}(V(h(x)))$, where 
\par $(vii)$ $\hat V_h=V\circ h$, 
$V$ is a vector field in $M$, $\hat V$ is a vector field 
in $Di^{ \{ l \} }_{ \{ \gamma \} , \delta ,\eta }(M).$
\par The proof in the complex case is analogous.
\par {\bf 2.9. Proposition.} {\it The loop group $G:=(L^MN)_{\xi }$
from \S 2.1.5 and the diffeomorphism groups $G:=Diff^t_{\beta ,\gamma }(M)$
from \cite{lurim1} and $G:=Di^{ \{ l \} }_{ \{ \gamma \} ,\delta ,\eta }(M)$ 
from \S 2.3 and $G:=Diff^{\xi }(M)$ from \S 2.2.1 have uniform atlases.}
\par {\bf Proof.} In view of Theorems 3.1 and 3.3 \cite{lurim1}
and Theorems 2.9.(1-4) \cite{lulgcm} and Lemma 2.6 and Theorem 2.8
above the diffeomorphism groups 
$G$ and the loop group $G:=(L^MN)_{\xi }$ have uniform atlases (see \S 2.1)
consistent with their topology, where $M$ is the real 
manifold $1\le t<\infty $,
$0\le \beta <\infty $, $0\le \gamma \le \infty $
for the diffeomorphism group $Diff^t_{\beta ,\gamma }(M)$
(see \cite{lurim1}). Others parameters are specified in the cited paaagraphs. 
They also include the particular cases of finite dimensional 
manifolds $M$ and $N$.
\par The case of complex compact $M$ for $G:=Diff^{\infty }(M)$ is trivial,
since $Diff^{\infty }(M)$ is the finite dimensional Lie group 
for such $M$ \cite{kobtg}. 
\par In view of Theorems 2.1.7, 2.8 
and Formulas 2.8.$(vi,vii)$ above and
Theorem 3.3 \cite{lurim1} that to satisfy conditions $(U1,U2)$
of \S 2.1.1 it is sufficient to find an atlas $At(G)$
of each such group $G$, for which $U_1$ is a neighbourhood of $e$,
$U_x^1$ and $U_x^2$ are for $x=e$ such that $\phi _1(U_1)$
contains a ball of radius $r>0$. Due to the existence of the 
left-invariant metrics in each such topological groups
and its paracompactness and separability 
we can take a locally finite covering
$\{ U_j: g_j^{-1}U_j\subset U_1 : j\in {\bf N} \} $, where
$\{ g_j: j\in {\bf N } \} $ is a countable subset of pairwise 
distinct elements of the group, $g_1=e$. Using
uniform continuity of $\tilde E$ we can satisfy $(U1,U2)$ 
with $r>0$, since the manifolds $M$ for diffeomorphism groups
and $N$ for loop groups also have uniform
atlases. Choosing $U_1$ in addition such that 
$\tilde E$ is bounded
on $U_1U_1$ and using left shifts $L_hg:=hg$, where $h$ and $g\in G$,
$AB:=\{ c: c=ab, a\in A, b\in B \} $ for $A\cup B\subset G$, 
and Condition $(U3)$ for $M$ and $N$ we get, that  
there exist sufficiently small neighbourhoods $U_1$, 
$U_e^1$ and $U_e^2$ with $U_e^2U_e^2\subset U_e^1$
and $U_x^1\subset xU_e^1$, $U_x^2\subset xU_e^2$ 
for each $x\in G$ such that Conditions $(U1-U3)$ are satisfied,
since uniform atlases exist on the Banach or Hilbert 
tangent space $T_eG$.
\section{Differentiable transition Wiener measures
on loop and diffeomorphism groups.}
\par {\bf 3.1. Definitions and Notes.}  Let $G$ be a 
Hausdorff topological group, we denote by 
$\mu : Af(G,\mu )\to [0,\infty )\subset \bf R$ 
a $\sigma $-additive measure. Its left shifts 
$\mu _{\phi }(E):=\mu (\phi ^{-1}\circ E)$ are considered for each 
$E \in Af(G,\mu )$,
where $Af(G,\mu )$ is
the completion of $Bf(G)$ by $\mu $-null sets, 
$Bf(G)$ is the Borel $\sigma $-field on $G$,
$\phi \circ E:=\{ \phi \circ h: h\in E \} $, $\phi \in G$.
For a monoid or a groupoid $G$ let left shifts of a measure
$\mu $ be defined by the following formula:
$\mu _{\phi }(E):=\mu (\phi \circ E)$.
Then $\mu $ is called quasi-invariant if there exists a dense subgroup
$G'$ (or submonoid or subgroupoid correspondingly) such that 
$\mu _{\phi }$ is equivalent to $\mu $ for each $\phi \in G'$.
Henceforth, we assume that a
quasi-invariance factor $\rho _{\mu }(\phi ,g)=\mu _{\phi }(dg)/\mu (dg)$
is continuous by $(\phi ,g) \in G' \times G$,
$\mu (V)>0$ for some (open) neighbourhood $V\subset
G$ of the unit element $e \in G$ and
$\mu (G)<\infty $. 
\par Let $({\sf M,F})$ be a space $\sf M$ of measures on $(G,Bf(G))$
with values
in $\bf R$ and $G"$ be a dense subgroup (or submonoid or
subgroupoid) in $G$ such that a topology $\sf F$ on
$\sf M$ is compatible with $G"$, that is, $\mu \mapsto \mu _h$
is the homomorphism
of $({\sf M,F})$ into itself for each $h \in G"$. Let $\sf F$ be the
topology of convergence for each $E \in Bf(G)$.
Suppose also that $G$ and $G"$ are real Banach manifolds such that
the tangent space $T_eG"$ is dense in $T_eG$, then $TG$ and $TG"$
are also Banach manifolds. Let $\Xi (G")$ denotes
the set of all differentiable vector fields $X$ on $G"$, that is, 
$X$ are sections of the tangent bundle $TG"$. We say that the measure
$\mu $ is continuously differentiable if there exists its tangent mapping
$T_{\phi }\mu _{\phi }(E)(X_{\phi })$ corresponding to the strong
differentiability relative to the Banach structures of the manifolds 
$G"$ and $TG"$. Its differential we denote $D_{\phi }\mu _{\phi }(E)$, 
so $D_{\phi }\mu _{\phi }(E)(X_{\phi })$ is the $\sigma $-additive
real measure by subsets $E\in Af(G, \mu )$ for each $\phi \in G"$ 
and $X\in \Xi (G")$ such that $D\mu (E): TG"\to \bf R$ is continuous
for each $E\in Af(G, \mu )$, 
where $D_{\phi }\mu _{\phi }(E)=pr_2\circ (T\mu )_{\phi }(E)$,
$pr_2: p\times {\bf F}\to \bf F$ is the projection in $TN$, $p\in N$,
$T_pN=\bf F$, $N$ is another real Banach differentiable manifold modelled on 
a Banach space $\bf F$, for a differentiable mapping $V: G"\to N$ 
by $TV: TG"\to TN$ is denoted the corresponding tangent mapping,
$(T\mu )_{\phi }(E):=T_{\phi }\mu _{\phi }(E)$. 
Then by induction $\mu $ is called $n$ times 
continuously differentiable if $T^{n-1}\mu $ is continuously
differentiable such that 
$T^n\mu :=T(T^{n-1}\mu )$, $(D^n\mu )_{\phi }(E)(X_{1, \phi },
...,X_{n, \phi })$ are the $\sigma $-additive real 
measures by $E\in Af(G, \mu )$
for each $X_1,$...,$X_n\in \Xi (G")$, where $(X_j)_{\phi }=:X_{j, \phi }$
for each $j=1,...,n$ and $\phi \in G"$, $D^n\mu : Af(G, \mu )\otimes
\Xi (G")^n\to \bf R$.
\par Differentiable quasi-invariant transition measures on
loop and diffeomorphism groups $G$ relative to dense subsgroups $G'$ are 
given by the following theorem, where the dense subgroups $G'$ are 
described precisely.
\par {\bf 3.2. Note.} Suppose that in the either
$Y^{\Upsilon ,b}$-Hilbert or $Y^{\Upsilon ,b, d'}$-manifold
$N$ modelled on
$l_2$ (see \S 2.1) there exists a dense 
$Y^{\Upsilon ,b'}$- or $Y^{\Upsilon ,b',d"}$-Hilbert 
submanifold $N'$ modelled on $l_{2,\epsilon }
=l_{2,\epsilon }({\bf C})$ (see \S 2.2.2), where 
\par $(1)$ $a>b>b'$ and $c>d'$ and either
\par $(2)$  $\infty >\epsilon >1/2$ and $d'\ge d"$  or
\par $(3)$ $\infty >\epsilon \ge 0$ and $d'>d"$
(such that either $d'_1>d"_1$ or $d'_1=d"_1$ and $d'_2>d"_2+1$)
correspondingly.
\par If $N$ is
finite dimensional let $N'=N$. Evidently, each $Y^{\Upsilon ,b}$-manifold
is the complex $C^{\infty }$-manifold. Certainly we suppose,
that a class of smoothness of a manifold $N'$ is not less than 
that of $N$ and classes of smoothness of $M$ and $N$ are not less 
than that of a given loop group for it
as in \S 2.1.5 and of $G'$ as below.
For the chosen loop group $G=(L^MN)_{\xi }$ let its dense subgroup 
$G':= (L^MN')_{\xi '}$ be the same
as in Theorem 2.11 \cite{lulgcm} or Theorem 2.6 \cite{lulgrm}
or \cite{ludan} with parameters:
\par $(a)$ $\xi '=(\Upsilon ,a")$ such that 
$a">b$
for $\xi ={\sf O}$ and the $Y^{\Upsilon ,b}$-manifolds $M$ and $N$
and the $Y^{\Upsilon ,b'}$-manifold $N'$;
\par $(b)$ $\xi '=(\Upsilon ,a")$ such that
$a>a">b$ for $\xi =(\Upsilon ,a)$; 
\par $(c)$ $\xi '=(\Upsilon ,a",c")$ for $\xi =(\Upsilon ,a,c)$
and $dim_{\bf C}M=\infty $
such that $b<a"<a$ and $d'<c"<c$ and either
$(2)$ $\infty >\epsilon >1$ with $d"\le d'$ or 
$(3)$ $\infty >\epsilon \ge 0$ with $d"<d'$,
such that either $d'_1>d"_1$ or $d'_1=d"_1$ and $d'_2>d"_2+1$,
where $M$ and $N$ are $Y^{\Upsilon ,b,d'}$-manifolds,
$N'$ is the $Y^{\Upsilon ,b',d"}$-manifold,
$1\le dim_{\bf C}M=:m<\infty $ in the cases $(a-b)$,
where either $a_1>a"_1$ or $a_1=a"_1$ with $a_2>a"_2+1$,
analogously for $c$ and $c"$, $b$ and $b'$ instead of $a$ and $a"$.
For the corresponding pair 
$G':=(L^M_{\bf R}N')_{\xi '}$ and $G:=(L^M_{\bf R}N)_{\xi }$
let indices in $(1-3)$ and $(a-c)$ be the same with substitution
of $\xi =\sf O$ on $\xi =(\infty ,H)$.
For real manifolds $M$ and $N$ in addition $N'$ is on 
$l_{2,\epsilon }({\bf R})$.
\par Then the embedding $J: T_eG'\hookrightarrow T_eG$
is the Hilbert-Schmidt operator, that follows
from \S \S 2.1 and 2.7.
\par For the diffeomorphism group $Diff^t_{\beta ,\gamma }(\tilde M)$
of a Banach manifold $\tilde M$ let $M$ be a dense Hilbert submanifold
in $\tilde M$ as in \cite{lurim1,lurim2}.
\par {\bf 3.3. Theorem.} {\it Let $G$ be either a loop group
or a diffeomorphism group for real or complex separable metrizable
$C^{\infty }$-manifolds $M$ and $N$,
then there exist a Wiener process on $G$ which induce
quasi-invariant infinite differentiable measures $\mu $
relative to a dense subgroup $G'$.}
\par {\bf Proof.} These topological groups also have structures 
of $C^{\infty }$-manifolds, but they do not satisfy the 
Campbell-Hausdorff formula in any open local subgroup.
Their manifold structures and actions of $G'$ on $G$ will be sufficient
for the construction of desired measures. Manifolds over $\bf C$ 
naturally have structures of manifolds over $\bf R$ also.
\par We take $G=\bar G$ and $Y=\bar Y$ for each loop group
$(L^MN)_{\xi }$ outlined in 3.2.$(b,c)$, for each diffeomorphism group 
$Diff^{\xi }(M)$ of a complex manifold $M$
given above, for each diffeomorphism group 
$G:=Di^{ \{ l \} }_{ \{ \gamma \} , \delta ,\eta }(M)$ 
for a real manifold $M$,
since such $G$ has the Hilbert manifold structure
(see Theorems 2.1.7 and 2.8 and also Appendix).
For ${\bar G}:=Diff^t_{\beta ,\gamma }(\tilde M)$ there exists
a Hilbert dense submanifold $M$ in a Banach manifold $\tilde M$ 
(see \S 2.6) and a subgroup
$G:=Di^{ \{ l \} }_{ \{ \gamma \} , \delta ,\eta }(M)$  
dense in $\bar G$ and a diffeomorphism subgroup
$G'$ dense in $G$ (see the proof of Theorem 3.10 \cite{lurim2}
and Lemma 2.6.2 above).
This $G'$ can be chosen as in Lemma 2.7. 
\par For the chosen loop group $G=(L^MN)_{\xi }$ let its dense subgroup 
$G':= (L^MN')_{\xi '}$ be the same
as in \S 3.2 Cases $(b,c)$. In case $3.2.(a)$
let $\bar G=(L^MN)_{\xi }$ and $G=(L^MN')_{\hat \xi }$
with $\hat \xi =(\Upsilon ,\hat a)$ such that $\hat a>a"$,
then $G'$ let be as in $3.2.(a)$.
\par On $G'$ there exists a $1$-parameter group 
$\rho : {\bf R}\times G'\to G'$
of diffeomorphisms of $G'$ generated by a $C^{\infty }$-vector field 
$X_{\rho }$ on $G'$ such that $X_{\rho }(p)=(d\rho (s,p)/ds )|_{s=0}$,
where $\rho (s+t,p)=\rho (s,\rho (t,p))$ for each $s, t\in \bf R$,
$\rho (0,p)=p$, $\rho (s,*): G'\to G'$ is the diffeomorphism
for each $s\in \bf R$ (about $\rho $ see \S 1.10.8 \cite{kling}).
Then each measure $\mu $ on $G$ and $\rho $ produce a $1$-parameter
family of measures $\mu _s(W):=\mu (\rho (-s,W))$.
For the construction of differentiable measures on the 
$C^{\infty }$-manifold we shall use the following statement:
if $a\in C^{\infty }(TG',TG)$ and $A\in C^{\infty }
(TG',L_{1,2}(TG',TG))$ and $a_x\in T_xG$ and $A_x\in L_{1,2}(T_xG',T_xG)$
for each $x\in G'$, each derivative by $x\in G'$:
$a^{(k)}_x$ and $A^{(k)}_x$ is a
Hilbert-Schmidt mappings into $Y=T_eG$ for each 
$k\in \bf N$ and $\sup_{\eta \in G} \| A_{\eta }(t)A_{\eta }^*(t) \| ^{-1}
\le C$, where $C>0$ is a constant,
then the transition probability $P(\tau ,x,t,W):=
P \{ \omega : \xi (t,\omega )=x, \xi (t,\omega )\in W \} $ 
is continuously stronlgy $C^{\infty }$-differentiable 
along vector fields on $G'$, where $G'$ is  
a dense $C^{\infty }$-submanifold on a space $Y'$,
where $Y'$ is a separable real Hilbert space having embedding into $Y$
as a dense linear subspace (see Theorem 3.3 and the Remark
after it in Chapter 4 \cite{beldal}
as  well as Theorems 4.2.1, 4.3.1 and 5.3.3 \cite{beldal}, 
Definitions 3.1 above), $W\in {\sf F}_t$.
\par Now let $G$ be a loop or a diffeomorphism group of 
the corresponding manifolds 
over the field $\bf R$ or $\bf C$. Then $G$ has the manifold structure.
If $exp^N: {\tilde T}N\to N$
is an exponential mapping of the manifold $N$, then
it induces the exponential $C^{\infty }$-mapping
${\tilde E}: {\tilde T}(L^MN)_{\xi }\to (L^MN)_{\xi }$
defined by ${\tilde E}_{\eta }(v)=exp^N_{\eta }\circ v_{\eta }$
(see Theorem 2.1.7),
where ${\tilde T}N$ is a neighbourhood of $N$ in a tangent
bundle $TN$, $\eta \in (L^MN)_{\xi }=:G$, $W_e$ is a neighbourhood 
of $e$ in $G$, $W_{\eta }=W_e\circ \eta .$
At first this mapping is defined for classes of equivalent 
mappings of the loop monoid $(S^MN)_{\xi }$ 
and then on elements of the group, since
$exp^N_{f(x)}$ is defined for each $x\in M$ and 
$f\in \eta \in (S^MN)_{\xi }$
(see Theorem 2.9.(3) \cite{lulgcm} and Theorem 2.4.3 \cite{lulgrm}).
The manifolds $G$ and $G'$ are of class $C^{\infty }$
and the exponential mappings $\tilde E$ and $\bar E$ for $G$ and $G'$
correspondingly are of class (strongly) $C^{\infty }$.
The analogous connection there exists in the diffeomorphism group
of the manifold $M$ satisfying the corresponding conditions
(see Theorem 3.3 \cite{lurim1}, \S 2.3 and Theorem 2.8)
for which: ${\tilde E}_{\eta }(v)=exp_{\eta (x)}\circ v_{\eta }$ for each 
$x\in M$ and $\eta \in G$. 
We can choose the uniform atlases $At_u(G)$ such that Christoffel symbols
$\Gamma _{\eta }$ are bounded on each chart
(see Proposition 2.9).
This mapping $\tilde E$ is for $G$ as the manifold and has not relations with
its group structure such as given by the Campbell-Hausdorff formula
for some Lie group, for example, finite dimensional Lie group. 
For the case of manifolds $M$ and $N$ over $\bf C$
we consider $G$ and others appearing manifolds with their structure 
over $\bf R$, since ${\bf C}={\bf R}\oplus i \bf R$
as the Banach space over $\bf R$.
\par Then for the manifold $G$ there exists an It$\hat o$ bundle.
Consider for $G$ an It$\hat o$ field $\sf U$ with a principal part
$(a_{\eta },A_{\eta })$, where $a_{\eta }\in T_{\eta }G$ and
$A_{\eta } \in L_{1,2}(H,T_{\eta }G)$ and $ker (A_{\eta })= \{ 0 \} $, 
$\theta : H_G\to G$ is a trivial bundle 
with a Hilbert fiber $H$ and $H_G:=G\times H$,
$L_{1,2}(\theta ,\tau _{\eta })$
is an operator bundle with a fibre $L_{1,2}(H,T_{\eta }G)$. 
To satisfy conditions of quasi-invariance and
differentiability of transition measures theorem 
we choose $A$ also such that $\sup_{\eta \in G} \| A_{\eta }(t)
A_{\eta }^*(t) \| ^{-1} \le C$, where $C>0$ is a constant.
If an operator $B$ is selfadjoint, then
$A_{\eta }^{\phi }B{A_{\eta }^{\phi }}^*$ is also selfadjoint, where
$A_{\eta }(t)=:A_{\eta }^{\phi _j}(t)$ is on a chart $(U_j,\phi _j)$.
If $\mu _B$ is a Gaussian measure on $T_{\eta }G$ with the correlation 
operator $B$, then $\mu _{A_{\eta }^{\phi }B{A_{\eta }^{\phi }}^*}$
is the Gaussian measure on $X_{1,\eta },$
where $B$ is selfadjoint and $ker (B)= \{ 0 \} $,
$A_{\eta }: T_{\eta }G\to X_{1,\eta }$, $X_{1,\eta }$ is a Hilbert space. 
We can take initially $\mu _B$ a cylindrical measure on
a Hilbert space $X'$ such that
$T_{\eta }G'\subset X'\subset T_{\eta }G$. 
If $A_{\eta }$ is the Hilbert-Schmidt operator
with $ker (A_{\eta })= \{ 0 \} $, then
$A_{\eta }^{\phi }B{A_{\eta }^{\phi }}^*$ is nondegenerate selfadjoint
linear operator of trace class and such the so called Radonifying 
operator $A_{\eta }^{\phi }$ gives 
the $\sigma $-additive measure 
$\mu _{A_{\eta }^{\phi }B{A_{\eta }^{\phi }}^*}$
in the completion $X'_{1,\eta }$ of $X'$
with respect to the norm $\| x \| _1:= \| A_{\eta }x \| $
(see \S II.2.4 \cite{dalfo}, \S I.1.1 \cite{sko}, 
\S II.2.4 \cite{oksb}). Then
using cylinder subsets we get a new Gaussian $\sigma $-additive
measure on $T_{\eta }G$, which we denote also by
$\mu _{A_{\eta }^{\phi }B{A_{\eta }^{\phi }}^*}$
(see also Theorems I.6.1 and III.1.1 \cite{kuo}).
\par If $U_j\cap U_l\ne \emptyset $, then
$A_{\eta }^{\phi _l}(t)={f_{\phi _l,\phi _j}}'
A_{\eta }^{\phi _l}(t){{f_{\phi _l,\phi _j}}'}^{-1}$, hence
the correlation operator $A_{\eta }^{\phi }B{A_{\eta }^{\phi }}^*$
is selfadjoint on each chart of $G$, that produces the Wiener 
process correctly.
Therefore, we can consider a stochastic process
$\mbox{ }d\xi (t,\omega )={\tilde E}_{\xi (t,\omega )}
[a_{\xi (t,\omega )}dt+A_{\xi (t,\omega )}dw],$
where $w$ is a Wiener process on $T_{\eta }G$
defined with the help of nuclear nondegenerate 
selfadjoint positive definite operator $B$.
The corresponding Gaussian measures 
$\mu _{tA_{\eta }^{\phi }B{A_{\eta }^{\phi }}^*}$ for $t>0$
(for the Wiener process)  are defined on the Borel
$\sigma $-algebra of $T_{\eta }G$ and 
$\mu _{tA_{\eta }^{\phi }B{A_{\eta }^{\phi }}^*}$
for such Hilbert-Schmidt nondegenerate linear operators
$A_{\eta }$ with $ker (A_{\eta })= \{ 0 \} $
are $\sigma $-additive (see Theorem II.2.1 \cite{dalfo}).
When the embedding operator $T_{\eta }G'\hookrightarrow T_{\eta }G$
is of Hilbert-Schmidt class, then there exists $A_{\eta }$ and $B$
such that $\mu _{tA_{\eta }^{\phi }B{A_{\eta }^{\phi }}^*}$
is the quasi-invariant and $C^{\infty }$-differentiable
measure on $T_{\eta }G$
relative to shifts on vectors from $T_{\eta }G'$
(see Theorem 26.2 \cite{sko} using Carleman-Fredholm determinant
and Chapter IV \cite{dalfo} and \S 5.3 \cite{ustzak}).
Henceforth we impose such 
demand on $B$ and $A_{\eta }$ for each $\eta \in G'$.
\par Consider left shifts $L_h: G\to G$
such that $L_h\eta :=h\circ \eta $. 
Let us take $a_e\in T_eG$, $A_e\in L_{1,2}(T_eG',T_eG)$, 
then we put $a_x=(DL_x)a_e$ and $A_x=(DL_x)\circ A_e$
for each $x\in G$, hence $a_x\in T_eG$ and 
$A_x\in L_{1,2}(H_x,(DL_x)T_eG)$, where $(DL_x)T_eG=T_xG$
and $T_eG'\subset T_eG$, $H_x:=(DL_x)T_eG'$.
Operators $L_h$ are (strongly) $C^{\infty }$-differentiable
diffeomorphisms of $G$ such that $D_hL_h: T_{\eta }G\to 
T_{h\eta }G$ is correctly defined, since $D_hL_h=h_*$
is the differential of $h$ \cite{ebi,eichh}.
In view of the choice of $G'$ in $G$ each 
covariant derivative $\nabla _{X_1}...\nabla _{X_n}(D_hL_h)Y$
is of class $L_{n+2,2}({TG'}^{n+1}\times G',TG)$ 
for each vector fields $X_1,...,X_n,Y$ on $G'$
and  $h\in G'$, since for each $0\le l\in \bf Z$ the embedding of $T^lG'$ 
into $T^lG$ is of Hilbert-Schmidt 
class, where $T^0G:=G$ (above and in \cite{beldal}
mappings of trace and Hilbert-Schmidt classes were defined for
linear mappings on Banach and Hilbert spaces and then 
for mappings on vector bundles). Take a dense subgroup $G'$
as it was otlined above and consider left shifts $L_h$ for $h\in G'$.
\par The considered here groups $G$ are separable, 
hence the minimal $\sigma $-algebra generated by cylindrical 
subalgebras $f^{-1}({\sf B}_n)$, n=1,2,..., coincides with 
the $\sigma $-algebra $\sf B$ of Borel subsets of $G$, where 
$f: G\to \bf R^n$ are continuous functions, ${\sf B}_n$ 
is the Borel $\sigma $-algebra  of $\bf R^n$. Moreover, $G$ 
is the topological Radon space (see Theorem I.1.2 and 
Proposition I.1.7 \cite{dalfo}).
Let $P(t_0,\psi ,t,W):=P( \{ \omega : \xi (t_0,\omega )=\psi ,
\xi (t,\omega )\in W \} )$ be the transition probability of 
the stochastic process $\xi $ for $0\le t_0<t$, which is defined on a 
$\sigma $-algebra $\sf B$ of Borel subsets in $G$, $W\in \sf B$,
since each measure
$\mu _{A_{\eta }^{\phi }B{A_{\eta }^{\phi }}^*}$
is defined on the $\sigma $-algebra
of Borel subsets of $T_{\eta }G$ (see above). 
On the other hand, $S(t,\tau ;gx)=gS(t,\tau ;x)$
is the stochastic evolution family of operators for each
$0\le t_0\le \tau <t$.
There exists $\mu (W):=P(t_0,\psi , t,W)$ such that 
it is a $\sigma $-additive
quasi-invariant strongly $C^{\infty }$-differentiable relative to 
the action of $G'$ by the left shifts $L_h$ on $\mu $
measure on $G$, for example, $t_0=0$ and $\psi =e$
with $t_0<t$, that is, $\mu _h(W):=\mu (h^{-1}W)$
is equivalent to $\mu $ and it is strongly infinitely differentiable 
by $h\in G'$.
\par The proof in cases $G=\bar G$ is thus obtained.
In cases $G\subset \bar G$ and $G\ne \bar G$ the use of the standard 
procedure of cylinder subsets induce a Weiener process and a transition 
measure from $G$ on $\bar G$ which is quasi-invariant and 
$C^{\infty }$-differentiable relative to $G'$ (see aslo 
\cite{lurim2}).
\par {\bf 3.4. Note.} This proof also shows, 
that $\mu $ is infinitely differentiable 
relative to each 1-parameter group $\rho : {\bf R}\times G'\to G'$
of diffeomorphisms of $G'$ generated by a $C^{\infty }$-vector field 
$X_{\rho }$ on $G'$. 
Evidently, considering different $(a,A)$ we see that there exist
${\sf c}=card ({\bf R})$ nonequivalent Wiener processes
on $G$ and $\sf c$ orthogonal quasi-invariant $C^{\infty }$-differentiable
measures on $G$ relative to $G'$ (see the Kakutani theorem in \cite{dalfo}).
\section{Differentiable Wiener transition measures
on loop monoids.}
\par This section is the consequence of the preceding 
sections and contains results for loop monoids
as well as for loop groupoids, which are defined 
in \S 4.2. For the considered here classes of manifolds
the generalized path space is defined in \S 4.4.
Differentiable transition Wiener measures on them
are given in Theorems 4.1, 4.3 and 4.5.
\par {\bf 4.1. Theorem.} {\it Let $G:=(S^MN)_{\xi }$ be a loop monoid
for both real or complex manifolds $M$ and $N$.
Then there exists a dense submonoid $G':=(S^MN')_{\xi '}$ 
and a stochastic process,
which generates quasi-invariant strongly $C^{\infty }$-differentiable
measure $\mu $ on $G$ relative to $G'$.}
\par The {proof} is quite analogous to that of
Theorem 3.3 with the help of definiton 3.1.
Pairs $(\xi ,\xi ')$ and $(N,N')$ were given above
in \S 3.2.
\par {\bf 4.2. Note and definition.} 
Let now $M$ and $N$ be two orientable Riemann manifolds
finite or infinite dimensional. If $M_m$ is a compact manifold and 
$f_{n,m}\in Y^{\xi }(M_m,N)$
has a rank $rank (f_{n,m}(x))=dim_{\bf R}T_xM_m$ for each
$x\in M_m$, then $f_{n,m}(M_m)$ is the $Y^{\xi }$-submanifold
in $N$ and on $f_{n,m}(M_m)$ there exists the Levi-Civit\`a connection
and the Riemann volume element $\nu _{n,m}$ as in \S 2.1.5.1 
such that $\nu _{n,m}(f_{n,m}(M_m))=1$. 
This induces local normal coordinates in $f_{n,m}(M_m)$. 
In particular, if $M_m=S^1$ 
we get the natural parameter corresponding to the 
length of an acr in a curve,
analogously in the multi-dimensional case. For each function
$f\in Y^{\xi }(M,N)$ there exists a sequence
$f_{n,m(n)}|_{M_m}\in Y^{\xi } (M_m,N)$ converging to $f$,
hence there are the natural coordinates for $f$, which are
mappings $\psi _f\in Y^{\xi }(B,N)$ and $h_f\in Y^{\xi }
(f(M),N)$ with $h_f\circ \psi _f=f$, where $B$ is the unit sphere
in $\bf R^m$ or $l_2$ over $\bf R$ correspondingly.
There exists an embedding $\xi ^*: Y^{\xi }(M\vee M,N)
\hookrightarrow Y^{\xi }(M,N)$ (see \cite{ludan,lulgcm}
and \S \S 2.1.4, 2.1.5 above).
In combination with the choice of the natural coordinates
we get the following continuous composition
$g\circ f$ in $G:=Y^{\xi }(M,s_0;N,y_0)$
such that $g\circ w_0=g$, that supplies $G$ with the groupoid structure
with the unity.
Let $G':=Y^{\xi '}(M,s_0;N',y_0)$ with $\xi '=(\Upsilon ,a",c")$
for $\xi =(\Upsilon ,a,c)$, where $b'<a'<a"<a$ and $d"<c'<c"<c$, 
$N'$ is a $Y^{\Upsilon ,b',d"}$-submanifold dense in $N$
(see also Conditions $(b,c)$ in \S 3.2).
Such space $Y^{\xi }(M,s_0;N,y_0)$ is called the generalized pinned 
loop space.
\par {\bf 4.3. Theorem.} {\it On the groupoid $G$ 
there exists a stochastic process generating a quasi-invariant
continuosly $C^{\infty }$-differentiable measure
$\mu $ relative to the dense subgroupoid $G'$ (see \S 4.2).}
\par {\bf Proof}. Since $N$ is the $C^{\infty }$-manifold,
then for each curve $f(t,x): {\bf R}\times M\to N$
of class $C^{\infty }$ by $t$ there exists $\partial ^l
f(t,x)/\partial t^l$ for each $l\in \bf N$, hence $T^lY^{\xi }(M,N)
=Y^{\xi }(M,T^lN)$ for each $l\in \bf N$
and $Y^{\xi }(M,N)$ is the $C^{\infty }$-manifold 
with the exponential mapping $(Exp^Y_gV)(x)=exp^N_{g(x)}\circ v(g(x))$ 
for each $x\in M$ (see Proposition 1.2.3 and Corollary 1.6.8 \cite{kling} 
and \cite{eliass}), where $V=v\circ g$ is the vector field
on $Y^{\xi }(M,N)$, $v$ is the vector field on $N$, 
$g\in Y^{\xi }(M,N)$. Therefore, $Exp^YV$ is of class $C^{\infty }$ by
Frech\'et on $\tilde TY^{\xi }(M,N)$. Then $Y^{\xi }(M,s_0;N,y_0)$ 
(see the notation in \S 2.1.5) is its closed 
$C^{\infty }$-submanifold with $g(s_0)=y_0$ and for it
the restriction $Exp^Y|_{\tilde TY^{\xi }(M,s_0;N,y_0)}$ 
also is defined and is of class $C^{\infty }$.
In view of \S 3.2 the embedding of $Y^{\xi '}(M,N')$
into $Y^{\xi }(M,N)$ is of Hilbert-Schmidt class.
Repeating almost the same arguments (without the use of $h^{-1}$)
for groupoids $G$ and $G'$ as in
Theorem 3.3 we get the proof of Theorem 4.3.
\par {\bf 4.4.} Let $Y^{\xi }(M,N)$ be as in \S 2.1.5, then
$Y^{\xi }(M,N)$ be called the generalized path space, where
a fixed mapping $\theta $ is omitted.
If $M_k=[0,1]^k$ are submanifolds in $M$, $k=1,2,...,$ such that
$\bigcup_kM_k$ is dense in $M$, then the subspace
$Y^{\xi }_l(M,N):=\{ f: f\in Y^{\xi }(M,N), f(x)=f(y)\mbox{ when }
x^k=y^k\quad (mod \quad 1)\mbox{ for each }k \} $
is called the loop space, where $x=(x^k: k=1,2,..., x^k\in {\bf R})
\in M$. Let $\xi $ and $\xi '$ be the same as in \S 4.2.
\par {\bf 4.5. Theorem.} {\it On $Y^{\xi }(M,N)$ and $Y^{\xi }_l(M,N)$ 
there exists a Wiener process such that it generates quasi-invariant 
measures relative to vector fileds of $Y^{\xi '}(M,N)$ and
$Y^{\xi '}_l(M,N)$, respectively.}
\par {\bf Proof.} $Y^{\xi }(M,N)$, $Y^{\xi }_l(M,N)$,
$Y^{\xi '}(M,N)$ and $Y^{\xi '}_l(M,N)$ are $C^{\infty }$-manifolds 
with of class $C^{\infty }$ exponential mappings, since
the exponential mapping $\tilde TY^{\xi }(M,N)$ generates
the corresponding restriction on $\tilde TY^{\xi }_l(M,N)$ also of 
class $C^{\infty }$ (see the proof in \S 4.3). 
They have uniform atlases. 
Here we can take $a\in TG$ and $A\in L_{1,2}(\theta ,\tau )$ 
(see also \S 3.3 without relations with $DL_h$).
Each vector field $X$ on $Y^{\xi '}=:G'$ generates the $1$-parameter
diffeomorphism group $\rho _X$ of $G'$ (see \S 3.4).
Then repeating the major parts of the proof of \S 3.3
without $L_h$ and so more simply, 
but using actions of vectors fields of $TG'$ by $\rho _X$ on $Y^{\xi }(M,N)$
or $Y^{\xi }_l(M,N)$ correspondingly we get the statement of this theorem, 
since $(D_X\rho _X)Y$ and $[(\nabla _X)^n(D_X\rho _X)]Y$
are of class $L_{n+2,2}((TG')^{n+2},TG)$ for each vector fields $X$ and $Y$ 
on $G'$ and each $n\in \bf N$, where $G:=Y^{\xi }$.
\section{Unitary representations associated with quasi-invariant measures.}
This section contains results for unitary representations
associated with quasi-invariant measures, which may be in particular 
transition Wiener measures. The generalization 5.1.2 of theorems
from preceding works \cite{lurim2,lulgcm} is proved.
It is applied in \S 5.2 
to the considered here case of Wiener transition measures.
Then applications to induced representations of nonlocally compact 
topological groups are given
having in mind the examples of constructed quasi-invariant measures
on loop groups and diffeomorphism groups.
\par {\bf 5.1.1. Note.} The transition measures $P=:\nu $ on $G$
induce strongly continuous unitary regular representations
of $G'$ given by the following formula:
$T_h^{\nu }f(g):=(\nu ^h(dg)/\nu (dg))^{1/2}f(h^{-1}g)$
for $f\in L^2(G,\nu ,{\bf C})=:H$, $T_h^{\nu }\in U(H)$,
$U(H)$ denotes the unitary group of the Hilbert space $H$. 
For the strong continuity of $T_h^{\nu }$
the continuity of the mapping $G'\ni h\mapsto \rho _{\nu }(h,g)\in 
L^1(G,\nu ,{\bf C})$ and that $\nu $ is the Borel measure
are sufficient, where $g\in G$, since $\nu $ is the 
Radon measure (see its definition in Chapter I \cite{dalfo}). 
On the other hand, the continuity of $\rho _{\nu }(h,g)
=\nu ^h(dg)/\nu (dg)$ by $h$
from the Polish group $G'$ into $L^1(G,\nu ,{\bf C})$ follows from
$\rho _{\nu }(h,g)\in L^1(G,\nu ,{\bf C})$  for each 
$h\in G'$ and that $G'$ is the topological subgroup of $G$.
In section 3 mostly Polish groups $\bar G$ and $G'$ were considered.
When $\bar G$ was not Polish it was used an embedding into $\bar G$
of a Polish subgroup $G$ such that $G'\subset G\subset \bar G$
and a measure on $G$ induces a measure on $\bar G$ with the help 
of an algebra of cylindrical subsets. So the considered 
cases of representations reduce to the case of Polish groups $(G',G)$.
\par More generally it is possible to consider
instead of the group $G$ a Polish topological space
$X$ on which $G'$ acts jointly continuously: $\phi : (G'\times X)\ni
(h,x)\mapsto hx=:\phi (h,x)\in X$, $\phi (e,x)=x$
for each $x\in X$, $\phi (v,\phi (h,x))=\phi (vh,x)$
for each $v$ and $h\in G'$ and each $x\in X$.
If $\phi $ is the Borel function, then
it is jointly continuous \cite{fidal}.
\par A representation $T: G'\to U(H)$ is called topologically irreducible,
if there is not any unitary operator (homeomorphism) $S$ on $H$ and a 
closed (Hilbert) subspace $H'$ in $H$ such that 
$H'$ is invariant relative to
$ST_hS^*$ for each $h\in G'$, that is, $ST_hS^*(H')\subset H'$.
\par A topological space $S$ is called dense in itself if
$S\subset S^d$, where $S^d$ is the derivative set of $S$,
that is, of all limit points $x\in cl(S\setminus \{ x \} ) $, $x\in S$, 
where $cl(A)$ denotes the closure of a subset $A$ in $S$
(see \S 1.3 \cite{eng}).
\par A measure $\nu $ on $X$ is called ergodic, if for
each $U\in Af(X,\nu )$ and $F\in Af(X,\nu )$ with $\nu (U)
\times \nu (F)\ne 0$
there exists $h\in G'$ such that $\nu ((h\circ E)\cap F)\ne 0$.  
\par {\bf 5.1.2. Theorem.} {\it Let $X$ be an infinite
Polish topological space 
with a $\sigma $-additive $\sigma $-finite nonnegative nonzero 
ergodic Borel measure
$\nu $ with  $supp(\nu )=X$ and
quasi-invariant relative to an infinite dense in itself 
Polish topological group
$G'$ acting on $X$ by the Borel function $\phi $. If 
\par $(i)$ $sp_{\bf C} \{\psi |
\quad \psi (g):=(\nu ^h(dg)/\nu (dg))^{1/2}, h\in G' \} $
is dense in $H$ and
\par $(ii)$ for each $f_{1,j}$ and $f_{2,j}$ in $H$, $j=1,...,n,$ 
$n\in \bf N$ and each $\epsilon >0$ there exists $h\in G'$ such that
$|(T_hf_{1,j},f_{2,j})| \le \epsilon |(f_{1,j},f_{2,j})|$,
when $|(f_{1,j},f_{2,j})|>0$.
Then the regular representation $T: G'\to U(H)$ is topologically
irreducible.}
\par {\bf Proof.}  From Condition $(i)$ it follows, that
the vector $f_0$ is cyclic, where $f_0\in H$ 
and $f_0(g)=1$ for each $g\in X$. 
In view of $card (X)\ge \aleph _0$ and the ergodicity of $\nu $
for each $n\in \bf N$ there are subsets $U_j\in Bf(X)$ and $g_j\in G'$
such that $\nu ((g_jU_j)\cap (\bigcup_{i=1,...,j-1,j+1,...,n}U_i))=0$
and $\prod_{j=1}^n\nu _j(U_j)>0$.
Together with Condition $(ii)$ this implies, 
that there is not any finite dimensional 
$G'$-invariant subspace $H'$ in $H$ such that
$T_hH'\subset H'$ for each $h\in G'$ and $H'\ne \{ 0 \}$.
Hence if there is a $G'$-invariant closed subspace $H'\ne 0$
in $H$ it is isomorphic with the subspace
$L^2(V,\nu ,{\bf C})$, where $V\in Bf(X)$ with $\nu (V)>0$. 
\par Let ${\sf A}_G$ denotes a $*$-subalgebra of 
an algebra ${\sf L}(H)$ of bounded linear operators on $H$
generated by the family of unitary operators 
$\{ T_h: h\in G' \} $. In view of the von Neumann
double commuter Theorem (see \S VI.24.2 \cite{fell})
${{\sf A}_G}"$ coincides with the weak and strong operator closures of
${\sf A}_G$ in ${\sf L}(H)$, where ${{\sf A}_G}'$
denotes the commuting algebra of ${\sf A}_G$ and ${{\sf A}_G}"=
({{\sf A}_G}')'$. 
\par Each Polish space is \v{C}ech-complete. By the 
Baire category theorem in a \v{C}ech-complete space $X$
the union $A=\bigcup_{i=1}^{\infty }A_i$ of a sequence of
nowhere dense subsets $A_i$ is a codense subset (see
Theorem 3.9.3 \cite{eng}). On the other hand, in view of Theorem
5.8 \cite{hew} a subgroup of a topological group
is discrete if and only if  it contains an isolated point.
Therefore, we can choose
\par $(i)$ a probability 
Radon measure $\lambda $ on $G'$ such that $\lambda $ has not any atoms and
$supp (\lambda )=G'$.
In view of the strong continuity of
the regular representation there exists the S. Bochner integral
$\int_XT_hf(g)\nu (dg)$ for each $f\in H$, which implies its existence 
in the weak (B. Pettis) sence. The measures $\nu $ and $\lambda $
are non-negative and bounded, hence $H\subset L^1(X,\nu ,{\bf C})$
and $L^2(G',\lambda ,{\bf C})\subset L^1(G',\lambda ,{\bf C})$
due to the Cauchy inequality. Therefore, we can apply below 
the Fubini Theorem (see \S II.16.3 \cite{fell}).
Let $f\in H$, then there exists a countable orthonormal base
$\{ f^j: j\in {\bf N} \} $ in $H\ominus {\bf C}f$. Then for each
$n\in \bf N$ the following set $B_n:=\{ q\in L^2(G',\lambda ,{\bf C} ):$
$(f^j,f)_H=\int_{G'}q(h)(f^j,T_hf_0)_H\lambda (dh)$ for $j=0,...,n \} $
is non-empty, since the vector $f_0$ is cyclic, where $f^0:=f$. 
There exists $\infty >R>\| f\|_H$ such that $B_n\cap B^R=:B^R_n$
is non-empty and weakly compact for each $n\in \bf N$, 
since $B^R$ is weakly compact, where
$B^R:=\{ q\in L^2(G',\lambda ,{\bf C} ): \| q\| \le R \} $
(see the Alaoglu-Bourbaki Theorem in \S (9.3.3) \cite{nari}).
Therefore, $B_n^R$ is a centered system of closed subsets
of $B^R$, that is, $\cap_{n=1}^mB^R_n\ne \emptyset $
for each $m\in \bf N$, hence it has a non-empty intersection, consequently,
there exists $q\in L^2(G',\lambda ,{\bf C})$ such that
$$(ii)\mbox{ }f(g)=\int_{G'}q(h)T_hf_0(g)\lambda (dh)$$ for $\nu $-a.e.
$g\in X$.
If $F\in L^{\infty }(X,\nu ,{\bf C})$, $f_1$ and $f_2\in H$,
then there exist $q_1$ and $q_2\in L^2(G',\lambda ,{\bf C})$
satisfying Equation $(ii)$. Therefore, 
$$(iii)\mbox{ }(f_1,Ff_2)_H=:c=
\int_X\int_{G'}\int_{G'}{\bar q}_1(h_1)q_2(h_2)\rho _{\nu }^{1/2}(h_1,g)
\rho _{\nu }^{1/2}(h_2,g)F(g)\lambda (dh_1)\lambda (dh_2)\nu (dg).$$
$$\mbox{Let }\xi (h):=\int_X\int_{G'}\int_{G'}{\bar q_1}(h_1)q_2(h_2)
\rho _{\nu }^{1/2}(h_1,g) \rho _{\nu }^{1/2}(hh_2,g)
\lambda (dh_1)\lambda (dh_2) \nu (dg).$$
Then there exists $\beta (h)\in L^2(G',\lambda ,{\bf C})$
such that 
\par $(iv)$ $\int_{G'}\beta (h)\xi (h)\lambda (dh)=c$.\\
To prove this we consider two cases. If $c=0$ it is sufficient
to take $\beta $ orthogonal to $\xi $ in $L^2(G',\lambda ,{\bf C})$. 
Each function $q\in L^2(G',\lambda ,{\bf C})$ 
can be written as $q=q^1-q^2+iq^3-iq^4$,
where $q^j(h)\ge 0$ for each $h\in G'$ and $j=1,...,4$,
hence we obtain the corresponding decomposition for $\xi $,
$\xi =\sum_{j,k}b^{j,k}\xi ^{j,k}$, where $\xi ^{j,k}$ corresponds to
$q_1^j$ and $q_2^k$, where $b^{j,k}\in \{ 1,-1,i,-i \}$. 
If $c\ne 0$ we can choose $(j_0,k_0)$ for which $\xi ^{j_0,k_0}\ne 0$
and 
\par $(v)$ $\beta $ is orthogonal to others $\xi ^{j,k}$ with 
$(j,k)\ne (j_0,k_0)$.\\ 
Otherwise, if $\xi ^{j,k}=0$ for each
$(j,k)$, then $q_l^j(h)=0$ for each $(l,j)$ and $\lambda $-a.e. $h\in G'$,
since 
$$\xi (0)=\int_X\nu (dg)(\int_{G'}{\bar q_1}(h_1)\rho _{\nu }^{1/2}
(h_1,g)\lambda (dh_1))(\int _{G'}q_2(h_2)
\rho _{\nu }^{1/2}(h_2,g)\lambda (dh_2))=0$$ 
and this implies $c=0$, which 
is the contradiction with the assumption $c\ne 0$.
Hence there exists $\beta $ satisfying conditions $(iv, v)$.
\par Let $a(x)\in L^{\infty }(X,\nu ,{\bf C})$, $f$ and $g\in H$, 
$\beta (h)\in L^2(G',\lambda ,{\bf C})$. Since $L^2(G',\lambda ,{\bf C})$ 
is infinite dimensional, then for each finite family of 
$a\in \{ a_1,...,a_m \} \subset L^{\infty }(X,\nu ,{\bf C})$,
$f\in \{ f_1,...,f_m \} \subset H$ there exists
$\beta (h)\in L^2(G',\lambda ,{\bf C})$, $h\in G'$, such that
$\beta $ is orthogonal to $\int_X{\bar f}_s(g)
[f_j(h^{-1}g)
(\rho _{\nu }(h,g))^{1/2}-f_j(g)]\nu (dg)$ for each $s,j=1,...,m$. Hence
each operator of multiplication on $a_j(g)$
belongs to ${{\sf A}_G}"$, since due to Formula $(iv)$
and Condition $(v)$ there exists $\beta (h)$ such that 
$$(f_s,a_jf_l)=\int_X\int_{G'}{\bar f}_s(g)\beta (h)(\rho _{\nu }
(h,g))^{1/2}f_l(h^{-1}g)\lambda (dh) \nu (dg)=$$
$$=\int_X\int_{G'} {\bar f}_s(g)
\beta (h)(T_hf_l(g))\lambda (dh)\nu (dg)\mbox{, }
\int_X{\bar f}_s(g)a_j(g)f_l(g)\nu (dg)=$$
$$=\int_X \int_{G'}{\bar f}_s(g)\beta (h)f_l(g)\lambda (dh)\nu (dg)=
(f_s,a_jf_l).$$
Hence ${{\sf A}_G}"$ contains 
subalgebra of all operators of multiplication on functions from
$L^{\infty }(X,\nu ,{\bf C})$.
With $G'$ and a Banach algebra $\sf A$ 
the trivial Banach bundle ${\sf B}={\sf A}\times G'$ is associative, in 
particular let ${\sf A}=\bf C$ (see \S VIII.2.7 \cite{fell}).
\par The regular representation $T$ of $G'$ gives rise to a canonical regular
$H$-projection-valued measure $\bar P$:
$\bar P(W)f=Ch_Wf$, where $f\in H$, $W\in Bf(X)$, $Ch_W$ 
is the characteristic function of $W$. Therefore, $T_h\bar P(W)=\bar P
(h\circ W)T_h$ for each $h\in G'$ and $W\in Bf(X)$, since
$\rho _{\nu }(h,h^{-1}\circ g)\rho _{\nu }(h,g)=1=\rho _{\nu }(e,g)$ 
for each $(h,g)\in G'\times X$, 
$Ch_W(h^{-1}\circ g)=Ch_{h\circ W}(g)$ and $T_h(\bar P(W)f(g))
=\rho _{\nu }(
h,g)^{1/2}\bar P(h\circ W)f(h^{-1}\circ g)$. Thus $<T,\bar P>$ is 
a system of imprimitivity for $G'$ over $X$, which is denoted 
${\sf T}^{\nu }$. This means that conditions
$SI(i-iii)$ are satisfied: 
\par $SI(i)$ $T$ is a unitary representation
of $G'$; 
\par $SI(ii)$ $\bar P$ is a regular 
$H$-projection-valued Borel measure on $X$ and 
\par $SI(iii)$ $T_h\bar P(W)=\bar P(h\circ W)T_h$ for all $h\in G'$ 
and $W\in Bf(X)$. 
\par For each $F\in L^{\infty }(X,\nu ,{\bf C})$ let $\bar \alpha _F$
be the operator in ${\sf L}(H)$ consisting
of multiplication by $F$: $\bar \alpha _F(f)=Ff$ for each $f\in H$. 
The map $F\to \bar \alpha _F$ is  an isometric $*$-isomorphism
of $L^{\infty }(X,\nu ,{\bf C})$ into ${\sf L}(H)$
(see \S VIII.19.2\cite{fell}). Therefore, Propositions 
VIII.19.2,5\cite{fell}
(using the approach of this particular case given above) are applicable
in our situation.
\par If $\bar p$ is a projection onto a closed ${\sf T}^{\nu }$-stable
subspace of $H$, then $\bar p$ commutes with all
$\bar P(W)$. Hence $\bar p$ commutes with multiplication by all
$F\in L^{\infty }(X,\nu ,{\bf C})$, so by \S VIII.19.2 \cite{fell}
$\bar p=\bar P(V)$, where $V\in Bf(X)$. Also $\bar p$ commutes with all
$T_h$, $h\in G'$, consequently, $(h\circ V)\setminus V$ and 
$(h^{-1}\circ V)\setminus V$ are $\nu $-null for each $h\in G'$, 
hence $\nu ((h\circ V)\bigtriangleup V)=0$ for all $h\in G'$. In view 
of ergodicity of $\nu $ and Proposition VIII.19.5 \cite{fell}
either $\nu (V)=0$ or $\nu (X\setminus V)=0$, hence
either $\bar p=0$ or $\bar p=I$, where $I$ is the unit operator.
Hence $T$ is the irreducible unitary representation.
\par {\bf 5.2. Theorem.} {\it On a loop or a diffeomorphism group $G$
there exists a stochastic process, which generates a quasi-invariant 
measure $\mu $ relative to a dense subgroup $G'$ such that the associated 
regular unitary representation $T^{\mu }: G' \to U(L^2(G,\mu ,{\bf C}))$ 
is irreducible.}
\par {\bf Proof.} From the construction of $G'$ and $\mu $ in Theorem 3.3
it follows that, if a function $f\in L^1(G,\mu ,{\bf C})$ 
satisfies the following condition
$f^h(g)=f(g)$ $(mod $ $\mu )$ by $g\in G$ for each $h\in G'$, 
then $f(x)=const $
$( mod $ $\mu )$, where $f^h(g):=f(hg)$, $g\in G$.  
\par Let $f(g)=Ch_U(g)$ be
the characteristic function of a subset $U$, $U\subset G$, $U\in Af(G,\mu )$,
then $f(hg)=1 $ $\Leftrightarrow g\in h^{-1}U$.  If $f^h(g)=f(g)$ is true by
$g\in G$ $\mu $-almost everywhere, then $\mu (\{ g\in G:  f^h(g)\ne f(g) \}
)=0$, that is $\mu ( (h^{-1}U)\bigtriangleup U)=0$, consequently, 
the measure $\mu $ 
satisfies the condition $(P)$ from \S VIII.19.5 \cite{fell}, where
$A\bigtriangleup B:=(A\setminus B)\cup (B\setminus A)$ 
for each $A, B\subset G$.
For each subset $E\subset G$ the outer measure is bounded,
$\mu ^*(E)\le 1$, since $\mu (G)=1$ and $\mu $ is non-negative, 
consequently, there exists $F\in
Bf(G)$ such that $F\supset E$ and $\mu (F)=\mu ^*(E)$.  
This $F$ may be interpreted as the representative
of the least upper bound in $Bf(G)$ relative to
the latter equality.
In view of the Proposition VIII.19.5 \cite{fell} the 
measure $\mu $ is ergodic.
\par In view of Theorems 2.1.7 and 2.8 the Wiener process on
the Hilbert manifold $G$ induces the Wiener process on 
the Hilbert space $T_eG$ with the help of the manifold
exponential mapping.
Then the left action $L_h$ of $G'$ on $G$
induces the local left action of $G'$ on a neighbourhood 
$V$ of $0$ in $T_eG$ with $\nu (V)>0$, where $\nu $ is induced by $\mu $.
A class of compact subsets
approximates from below each measure $\mu ^f$, 
$\mu ^f(dg):=|f(g)|\mu (dg)$,
where $f\in L^2(G,\mu ,{\bf C})=:H$.
Due to the Egorov Theorem II.1.11 \cite{fell} for each $\epsilon >0$
and for each sequence
$f_n(g)$ converging to $f(g)$ for $\mu $-almost every $g\in G$,
when $n\to \infty $, there exists a compact subset $\sf K$
in $G$ such that $\mu (G\setminus {\sf K})<\epsilon $ and
$f_n(g)$ converges on $\sf K$ uniformly by $g\in \sf K$,
when $n\to \infty $.
\par In view of Lemma IV.4.8 \cite{oksb}
the set of random variables $ \{ \phi (B_{t_1},...,B_{t_n}) :\quad
t_i \in [t_0,T], \phi \in C^{\infty }_0({\bf R^n}), n\in {\bf N} \} $
is dense in $L^2({\sf F}_T,\mu )$, where $T>t_0$. In accordance with
Lemma IV.4.9 \cite{oksb} the linear span of random variables
of the type $\{ exp \{ \int_0^Th(t)dB_t(\omega )-\int_0^Th^2(t)dt/2 \} :\quad
h\in L^2[t_0,T]$ (deterministic) $ \} $ is dense in $L^2({\sf F}_T,\mu )$,
where $T>t_0$. 
Therefore, in view of Girsanov Theorem 2.1.1 and Theorem 5.4.2 \cite{ustzak}
the following space 
$sp_{\bf C}\{ \psi (g):=(\rho (h,g))^{1/2}: h\in G' \}=:Q$ is dense in
$H$, since $\rho _{\mu }(e,g)=1$ for each $g\in G$
and $L_h: G\to G$ are diffeomorphisms of the manifold $G$, $L_h(g)=hg$.
Finally we get from Theorem 3.3 above that there exists
$\mu $, which is ergodic and Conditions $(i,ii)$ of Theorem 5.1.2
are satisfied. Evidently $G'$ and $G$ are infinite and dense in themselves.
Hence from Theorem 5.1.2 the statement of this theorem, follows.
\par {\bf 5.3. Note.} Then analogously to \S 3.3 there can be constructed 
quasi-invariant and pseudo-differentiable measures on the manifold
$M$ relative to the action of the diffeomorphism group 
$G_M$ such that $G'\subset G_M$. Then Poisson measures on 
configuration spaces associated with either $G$ or $M$ can be
constructed and producing new unitary representations 
including irreducible as in \cite{lupom}.
\par Having a restriction of a transition measure $\mu $ from \S 3.3
on a proper open neighbourhood of $e$ in $G$ it is possible
to construct a quasi-invariant $\sigma $-finite nonnegative
measure $m$ on $G$ such that $m(G)=\infty $ using left shifts
$L_h$ on the paracompact
$G$. Analogously such measure can be constructed on the 
manifold $M$ in the case of the diffeomorphism group using Wiener 
processes on $M$. For definite $\mu $ in view of Theorems 2.9 \cite{lupom}
and 5.2 the corresponding Poisson measure $P_m$ is ergodic.
Therefore, Theorems 3.4, 3.6, 3.9, 3.10, 3.13 and 3.14
\cite{lupom}
also encompass the corresponding class of measures $m$ and $P_m$
arising from the constructed in \S 3.3 transition measures.
\par In view of Proposition II.1 \cite{neeb} for the separable 
Hilbert space $H$ the unitary group endowed with the strong 
operator topology $U(H)_s$ is the Polish group.
Let $U(H)_n$ be the unitary group with the metric
induced by the operator norm. In view of the Pickrell's theorem 
(see \S II.2 \cite{neeb}): if $\pi : U(H)_n\to U(V)_s$
is a continuous representation of $U(H)_n$ on  
the separable Hilbert space $V$, then $\pi $ is also
continuous as a homomorphism from $U(H)_s$ into
$U(V)_s$. Therefore, if $T: G'\to U(H)_s$ is
a continous representation, then there are new representations
$\pi \circ T: G'\to U(V)_s$. On the other hand, the 
unitary representation theory of $U(H)_n$ is the same as that of
$U_{\infty }(H):=U(H)\cap (1+L_0(H))$, since the group $U_{\infty }(H)$
is dense in $U(H)_s$.
\par {\bf 5.4. Remark.} Let $\mu $ be a Borel regular Radon
non-negative quasi-invariant measure on a topological Hausdorff group $G$ 
relative to a dense subgroup $G'$ with a continuous quasi-invariance factor
$\rho _{\mu }(x,y)$ on $G'\times G$ and
$0<\mu (G)<\infty $.  
Suppose  that $V: S\to U(H_V)$ is a strongly continuous unitary 
representation of a closed subgroup $S$ in $G'$.
There exists a Hilbert space $L^2(G,\mu ,H_V)$ of equivalence classes 
of measurable functions $f: G\to H_V$ with a finite norm
$$(1)\mbox{ }\| f \| :=(\int_G \| f(g) \|^2_{H_V}\mu (dg))^{1/2}<\infty .$$
Then there exists a subspace $\Psi _0$ of functions
$f\in L^2(G,\mu ,H_V)$ such that $f(hy)=V_{h^{-1}}f(y)$ for each
$y\in G$ and $h\in S$, the closure of $\Psi _0$ in
$L^2(G,\mu ,H_V)$ is denoted by $\Psi ^{V,\mu }$. For each 
$f\in \Psi ^{V,\mu }$
there is defined a function
$$(2)\mbox{ }(T^{V,\mu }_xf)(y):=\rho _{\mu  }^{1/2}(x,y) f(x^{-1}y),$$
where $\rho _{\mu }(x,y):=\mu _x(dy)/\mu (dy)$ is a quasi-invariance 
factor for each $x\in G'$ and $y\in G$, $\mu _x(A):=
\mu (x^{-1}A)$ for each Borel subset $A$ in $G$.
Since $(T^{V,\mu }_xf)(hy)=V_{h^{-1}}((T_xf)(y))$, then $\Psi  ^{V,\mu }$
is a $T^{V,\mu }$-stable subspace. Therefore, $T^{V,\mu }: G'\to U(\Psi 
^{V,\mu })$
is a strongly continuous unitary representation,
which is called induced and denoted by $Ind_{S\uparrow G'}(V)$.
\par {\bf 5.5.1. Note.} Let $G$ be a topological Hausdorff group
with a non-negative
quasi-invariant measure $\mu $ relative to a dense subgroup $G'$.
Suppose that there are two closed
subgroups $K$ and $N$ in $G$ such that $K':=K\cap G'$ and $N'=N\cap G'$
are dense subgroups in $K$ and $N$ respectively. We say that
$K$ and $N$ act regularly in $G$, if there exists a sequence 
$\{ Z_i:$ $i=0,1,... \} $ of Borel subsets $Z_i$ satisfying two conditions:
\par $(i)$ $\mu (Z_0)=0$, $Z_i(k,n)=Z_i$ for each pair
$(k,n)\in K\times N$ and each $i$;
\par $(ii)$ if an orbit $\sf O$ relative to the action of
$K\times N$ is not a subset of $Z_0$, then
${\sf O}=\bigcap_{Z_i\supset \sf O}Z_i$, where $g(k,n):=k^{-1}gn$.
Let $T^{V,\mu }$ be a representation of
$G'$ induced by a unitary representation $V$ of $K'$ and
a quasi-invariant measure $\mu $ 
(for example, as in \S 3). 
We denote by $T^{V,\mu }_{N'}$ a restriction of $T^{V,\mu }$ on $N'$
and by $\sf D$ a space $K\setminus G/N$ of double coset classes $KgN$.
\par {\bf 5.5.2. Theorem.} {\it There are a unitary operator $A$ on 
$\Psi ^{V,\mu }$
and a measure $\nu $ on a space $\sf D$ such that 
$$(1)\mbox{ }A^{-1}T^{V,\mu }_nA=\int_{\sf D}T_n(\xi )d\nu (\xi )$$
for each $n\in N'$.
$(2).$ Each representation $N'\ni n\mapsto T_n(\xi )$ in 
the direct integral decomposition $(1)$ is defined relative to 
the equivalence of a double coset class $\xi $. For a subgroup
$N'\cap g^{-1}K'g$ its representations 
$\gamma \mapsto V_{g\gamma g^{-1}}$ are equivalent for each $g\in G'$
and representations $T_{N'}(\xi )$ can be taken up to their equivalence
as induced by $\gamma \mapsto V_{g\gamma g^{-1}}$.}
\par {\bf Proof.} A quotient mapping 
$\pi _{\sf X}: G\to G/K=:\sf X$
induces a measure $\hat \mu $ on $\sf X$ such that
${\hat \mu }(E)=\mu (\pi ^{-1}_{\sf X}(E))=:(\pi _{\sf X}^*\mu )(E)$ 
for each Borel subset 
$E$ in $\sf X$. In view of Radon-Nikodym theorem II.7.8 \cite{fell}
for each $\xi \in \sf D$
there exists a measure $\mu _{\xi }$ on $\sf X$ such that
$$(3)\mbox{ }d{\hat \mu }(x)= d\nu (\xi )d\mu _{\xi }(x),$$ 
where $x\in \sf X$, $\nu (E):=(s^*\mu )(E)$
for each Borel subset $E$ in $\sf D$, $s: G\to \sf D$ is
a quotient mapping.
In view of \S 26 \cite{nai} and Formula $(3)$
the Hilbert space $H^V:=L^2({\sf X},{\hat \mu },H_V)$ has a decomposition 
into a direct integral 
$$(4)\mbox{ } H^V=\int_{\sf D}H(\xi )d\nu (\xi ),$$
where $H_V$ denotes a complex Hilbert space of the representation
$V: K'\to U(H_V)$. Therefore, 
$$ \| f \| ^2_{H^V}=\int_{\sf D}\| f \| ^2_{H(\xi )}d\nu (\xi ).$$
From Formulas $(4)$ and $5.4.(1,2)$ we get the first statement 
of this theorem for a subspace $\Psi ^{V,\mu }$ of $H^V$.
\par If $f\in L^2({\sf X},{\hat \mu },H_V)$, then
$\pi ^*_{\sf X}f:=f\circ \pi _{\sf X}\in L^2(G,\mu ,H_V)$.
This induces an embedding $\pi ^*_{\sf X}$ of $H^V$ into
$\Psi ^{V,\mu }$. Let $\sf F$ be a filterbase of neighbourhoods $A$ 
of $K$ in $G$ such that $A=\pi ^{-1}_{\sf X}(S)$, where
$S$ is open in $\sf X$, hence $0<\mu (A)\le \mu (G)$
due to quasi-invaraince of $\mu $ on $G$ relative to $G'$.
Let $\psi \in \xi \in \sf D$, then $\psi =Kg_{\xi }$, where
$g_{\xi }\in G$, hence $\psi =\psi (N\cap g_{\xi }^{-1}Kg_{\xi })$.
In view of Formula $(3)$ for each $x\in N'$ and $\eta
=Kx$ we get
$$(5)\quad \rho _{\mu _{\xi }}^{1/2}(\eta ,\xi )=
\lim_{\sf F}[\int_A \rho ^{1/2}(x,zg_{\xi })d\mu (z)/\mu (A)],$$
since by the supposition $\rho _{\mu }(h,y)$ is continuous on 
$G'\times G$ (see also \S 1.6 \cite{eng} and \S \S 3.3, 5.1.1).
Therefore, 
$$(a,T_x(\xi )b)_{H^V}=\lim_{\sf F} [\int_A (\pi ^*a,
 \rho _{\mu }^{1/2}(x,zg_{\xi })(\pi ^*b)_x^{zg_{\xi }})_{
\Psi ^{V,\mu }}d\mu (z)/\mu (A)]$$
for each $x\in N'$ and $a,b \in H^V$, where $f_z^h(\zeta ):=
f(z^{-1}h\zeta )$ for a function $f$  on $G$ and $h, z, \zeta \in G$. 
In view of the cocycle condition $\rho _{\mu }(yx,z)=\rho _{\mu }
(x,y^{-1}z)\rho _{\mu }(y,z)$ for each $x, y \in G'$ and $z\in G$ we get 
$T_{yx}(\xi )=T_y(\xi )T_x(\xi )$ for each $x, y\in N'$ and $T_x(\xi )$ 
are unitary representations of $N'$. Then 
$$(a,T_{yx}(\xi )b)_{H^V}=
\lim_{\sf F} [\int_A (\pi ^*a,V_{g_{\xi }yg_{\xi }^{-1}}[
\rho _{\mu }^{1/2}(x,zg_{\xi })(\pi ^*b)_x^{zg_{\xi }}])_{\Psi ^{V,\mu }}
d\mu (z)/\mu (A)]$$ 
for each $y\in N'\cap g_{\xi }^{-1}K'g_{\xi }$.
Hence the representation $T_x(\xi )$ in the Hilbert space
$H(\xi )$ is induced by a representation 
$(N'\cap g_{\xi }^{-1}K'g_{\xi })\ni y\mapsto 
V_{g_{\xi }yg_{\xi }^{-1}}.$
\par {\bf 5.6.1. Note.} Let $V$ and $W$ be two unitary representations of 
$K'$ and $N'$ (see \S 5.5.1). In addition let $K$ and $N$ be regularly 
related in $G$ and $V{\hat \otimes }W$ denotes an external tensor product of
representations for a direct product group $K\times N$.
In view of \S 5.4 a representation $T^{V,\mu }{\hat \otimes }T^{W,\mu }$
of an external product group ${\sf G}:=G\times G$ 
is equvalent with an induced representation
$T^{V{\hat \otimes }W,\mu \otimes \mu }$, where $\mu \otimes \mu $
is a product measure on $\sf G$.
A restriction of $T^{V{\hat \otimes }W,\mu \otimes \mu }$
on ${\tilde G}:=\{ (g,g):$ $g\in G \}$
is equivalent with an internal tensor product
$T^{V,\mu }\otimes T^{W,\mu }$.
\par {\bf 5.6.2. Theorem.} {\it There exists a unitary operator $A$ 
on $\Psi ^{V{\hat \otimes }W,\mu \otimes \mu }$ such that
$$(1)\mbox{ }A^{-1}T^{V,\mu }\otimes T^{W,\mu }A=\int_{\sf D}
T(\xi )d\nu (\xi ),$$
where $\nu $ is an admissible measure on a space ${\sf D}:=
N\setminus G/K$ of double cosets.
\par $(2).$ Each representation $G'\ni g\mapsto T_g(\xi )$
in Formula $(1)$ is defined up to the equivalence of $\xi $ in $\sf D$.
Moreover, $T(\xi )$ is unitarily equivalent with
$T^{{\tilde V}\otimes {\tilde W},\mu \otimes \mu }$, 
where $\tilde V$ and $\tilde W$
are restrictions of the corresponding representations 
$y\mapsto V_{gyg^{-1}}$ and $y\mapsto W_{\gamma y\gamma ^{-1}}$ 
on $g^{-1}K'g\cap \gamma ^{-1}N'\gamma $,
$g, \gamma \in G'$, $g\gamma ^{-1}\in \xi $.}
\par {\bf Proof.} In view of \S 18.2 \cite{barut}
$P\setminus {\sf G}/{\tilde G}$ and $K\setminus G/N$ 
are isomorphic Borel spaces,
where $P=K\times N$. In view of Theorem 5.5.2
there exists a unitary operator $A$ on a subspace $\Psi ^{V
\hat \otimes W,\mu \otimes \mu }$ of the Hilbert space
$L^2({\sf G},\mu \otimes \mu ,H_V\otimes H_W)$ such that 
$$A^{-1}T^{V\hat \otimes W,\mu \otimes \mu }|_{\tilde G}A=
\int_{\sf D}T_{\tilde G}(\xi )d\nu (\xi ),$$
where each $T_{\tilde G}(\xi )$ is induced by
a representation $(y,y)\mapsto (V\hat \otimes W)_{(
g,\gamma )(y,y)(g,\gamma )^{-1}}$ of a subgroup
${\tilde G}'\cap (g,\gamma )^{-1}(K\times N)(g,\gamma )$,
the latter group is isomorphic with
$S:=g^{-1}K'g\cap \gamma ^{-1}N'\gamma ,$ that gives a representation
${\tilde V}{\hat \otimes }{\tilde W}$ of a subgroup $S\times S$ in $\sf G$.
Therefore, we get a representation $T^{{\tilde V}{\hat \otimes }{\hat W},
\mu \otimes \mu }$ equivalent with $Ind_{(S\times S)\uparrow {\sf G}'}
({\tilde V}{\hat \otimes }{\tilde W})|_{{\tilde G}'}$.
\par {\bf 5.7. Note.} Formulas $(3-5)$ from \S 5.5.2 also show how
a measure $\nu $ on a groupoid $Y^{\xi }(M,s_0;N,y_0)$ induces a measure
$\mu $ on $(S^MN)_{\xi }$ and produces an expression for a quasi-invariance 
factor on a loop monoid and then on a loop group.
\section{Appendix.}
Let us remind the principles
of the Wiener processes on manifolds.
\par Let $\bar G$ be a complete separable relative to its metric $\bar \rho $
$C^{\infty }$-manifold
on a Banach space $\bar Y$ over $\bf R$ such that it contains
a dense $C^{\infty }$-submanifold $G$ 
on a Hilbert space $Y$ over $\bf R$, where $G$ is also
separable and complete relative to its metric $\rho $.
Let $\tau _{G}: TG\to G$ be a tangent bundle on $G$. 
Let $\theta : Z_{G}\to G$ be a trivial bundle on $G$
with the fibre $Z$ such that $Z_{G}=Z\times G$, then $L_{1,2}(
\theta , \tau _{G})$ be an operator bundle with a fibre
$L_{1,2}(Z,Y)$, where $Z, Z_1,...,Z_n$ are Hilbert spaces, 
$L_{n,2}(Z_1,...,Z_n;Z)$ is
a subspace of a space of all Hilbert-Schmidt $n$ times multilinear 
operators from $Z_1\times ...\times Z_n$ into $Z$. Then 
$L_{n,2}(Z_1,...,Z_n;Z)$ has the structure of the Hilbert space
with the scalar product denoted by 
$$\sigma _2(\phi ,\psi ):=
\sum_{j_1,...,j_n=1}^{\infty }(\phi (e^{(1)}_{j_1},...,e^{(n)}_{j_n}),
\psi (e^{(1)}_{j_1},...,e^{(n)}_{j_n}))$$
for each pair of its elements $\phi , \psi .$ It
does not depend on a choise of the orthonormal bases
$\{ e{(k)}_j: j \} $ in $Z_k$. Let $\Pi :=
\tau _{G}\oplus L_{1,2}(\theta ,\tau _{G})$ be a Whitney sum
of bundles $\tau $ and $L_{1,2}(\theta ,\tau _{G})$. 
If $(U_j,\phi _j)$ 
and $(U_l,\phi _l)$ are two charts of $G$ with an open non-void intersection
$U_j\cap U_l$, then to a connecting mapping $f_{\phi _l,\phi _j}=
\phi _l\circ \phi _j^{-1}$ there corresponds a connecting mapping
$f_{\phi _l,\phi _j}\times {f'}_{\phi _l,\phi _j}$ for the bundle
$\Pi $ and its charts $U_j\times (Y\oplus L_{1,2}(Z,Y))$ for $j=1$ 
or $j=2$, where $f'$ denotes the strong derivative of $f$,
${f'}_{\phi _l,\phi _j}: (a^{\phi _j},A^{\phi _j})\mapsto
({f'}_{\phi _l,\phi _j}a^{\phi _j}, {f'}_{\phi _l,\phi _j}
\circ A^{\phi _j})$,
$a^{\phi }\in Y$ and $A^{\phi }\in L_{1,2}(Z,Y)$ 
for the chart $(U,\phi )$,
${f'}_{\phi _l,\phi _j}\circ A^{\phi _j}:=
{f'}_{\phi _l,\phi _j}A^{\phi _j}{f'}_{\phi _l,\phi _j}^{-1}$. 
Such bundles are called quadratic.
Then there exists a new bundle $J$ on $G$ with the same fibre as 
for $\Pi $, but with new connecting mappings: 
$J(f_{\phi _l,\phi _j}): (a^{\phi _j},A^{\phi _j})\mapsto
({f'}_{\phi _l,\phi _j}a^{\phi _j}+tr ({f"}_{\phi _l,\phi _j}(
A^{\phi _j}, A^{\phi _j}))/2, {f'}_{\phi _l,\phi _j}\circ A^{\phi _j})$,
where $tr (A)$ denotes a trace of an operator $A$.
Then using sheafs one gets the It$\hat o$ functor  $I: I(G)\to G$
from the category of manifolds to the category of quadratic bundles.
\par On a Hilbert space $W$ a 
distribution $\gamma _{b,B}$ is called Gaussian, 
if its Fourier transform is the following:
$$F'(\gamma _{b,B})(v)=exp \{ -(Bv,v)/2+i(b,v) \} ,$$
where $B$ is the corresponding symmetric
bounded nonnegative nondegenerate nuclear operator on $W$,
$b\in W$, $v\in W$.
\par On $Y$ let $B$ be a nuclear selfadjoint linear 
nonnegative operator with $ker (B)= \{ 0 \} $, 
then for each $t>0$ it defines 
a Gaussian measure $\mu _{tB}$ with zero mean and correlation
operator $tB$. It is defined with the help of
the Hilbert-Schmidt
structure in $Y$ (that is, the rigged Hilbert space):
$Y={Y'}_-,$ ${Y'}_+\hookrightarrow Y_0$ and
$Y_0\hookrightarrow {Y'}_-$ are Hilbert-Schmidt embeddings $B^{1/2}$, 
$Y':={Y'}_+$, where $(x,y)_+=(B^{-1}x,B^{-1}y)$ is the scalar product
in $Y'_+$ induced from the dense subspace $B^{-1}Y'_-$, 
$(x,y)_0=(B^{-1/2}x,B^{-1/2}y)$
is the scalar product in $Y_0$ induced from the dense subspace
$B^{-1/2}Y'_-,$ where $(x,y)$ is the scalar product in $Y=Y'_-$ 
for each $x, y \in Y$. By the definition a 
Wiener process $w(t,\omega )$ for $0\le t_0 \le t<\infty $
with values in $Y$ is a stochastic process for which
\par $(1)$ the differences $w(t_4,\omega )-w(t_3,\omega )$ 
and $w(t_2,\omega )-w(t_1,\omega )$ are independent
for each $t_0\le t_1<t_2\le t_3<t_4$; 
\par $(2)$ the random variable
$w(t+\tau ,\omega )-w(t,\omega )$ has a distribution 
$\mu _{\tau B}$, where $w(t_0,\omega ):=0$,
$(\Omega ,{\sf F},P)$
is a probability space of a set $\Omega $ with a $\sigma $-algebra $\sf F$
of its subset and a probability measure $P$. 
\par Then consider the class ${\sf K}(Y)$ of stochastic processes
$B(t,\omega )$ with values in $L_{1,2}({\sf H},Y)$ and satisfying
$\beta ^2(B)=\int_{t_0}^{\tau } {\sf M}\sigma _2(B(t, 
\omega ),B(t,\omega ))dt<\infty $, 
the space of all such operators is denoted by
$L_2(\Omega ,Y)$, where $\sf M$ denotes the operation 
of the mean value, the embedding of $\sf H$ into $Y$ is a 
Hilbert-Schmidt operator, $\omega \in \Omega $, 
$(\Omega ,{\sf F},P)$:
\par $(3)$ for each $t\ge t_0$ the quantity 
$B(t,\omega )$ is ${\sf F}_t$-measurable, where ${\sf F}_t$ 
is a flow of $\sigma $-algebras, that is, a monotone set of 
$\sigma $-algebras (${\sf F}_t\subset {\sf F}_s$ for each $s\ge t\ge t_0$)
such that for each $s\le t$ the random variable
$w(s,\omega )$ is ${\sf F}_t$-measurable, $w(\tau ,\omega )-w(s,\omega )$ 
is independent from ${\sf F}_t$ for each $\tau >s\ge t$.
\par Let ${\sf K}_0(Y)$ be the subset of ${\sf K}(Y)$ 
consisting of step functions $B(t,\omega )=B_j(\omega )$ for 
each $t_j\le t<t_{j+1}$, where $t_0<t_1<...<t_n=\tau $ is a 
partition of the segment $[t_0,\tau ]$ in $\bf R$.
In ${\sf K}_0(Y)$ the It$\hat o$ stochastic integral
is defined by ${\sf I}(B)=\int_{t_0}^{\tau }B(t,\omega )dt=\sum_{j=0}^{n-1}
B_j(\omega )[w(t_{j+1},\omega )-w(t_j,\omega )]$. It has the extension
${\sf I}: {\sf K}(Y)\to L_2(\Omega ,Y)$.
Let $a(t,\omega )$ be an ${\sf F}_t$-measurable function
with values in $Y$ such that $\int_{t_0}^{\tau }{\sf M} 
\| a(t,\omega )\|^2dt 
<\infty $ and let $\xi _0(\omega )$ be an ${\sf F}_{t_0}$-measurable
random variable. A stochastic process of the type
$$\xi (t,\omega )=\xi _0(\omega )+\int_{t_0}^t a(s,\omega )ds 
+ \int_{t_0}^t B(s,\omega )dw (s,\omega )$$
is said to have a stochastic differential and it is written as follows:
$$d\xi =a(t,\omega )dt+B(t,\omega )dw (t,\omega ).$$
If $f(t,x)$ is continuously differentaible by $t$ and  
continuously twice strongly differentiable by $x$ function
from $[t_0,\tau ]\times Y$ into $Y$ and they are bounded, then
$$f(t,\xi (t,\omega ))=f(t_0,\xi _0(\omega ))+
\int_{t_0}^t\{ {f'}_s(s,\xi (s,\omega ))+{f'}_x(s,\xi (s))a(s,\omega )+$$
$$tr (B^*(s,\omega ){f"}_{x,x}(s,\xi (s,\omega ))B(s,\omega ))/2 \} ds
+\int_{t_0}^t{f'}_x(s,\xi (s,\omega ))B(s,\omega )dw(s,\omega )$$
in accordance with the It$\hat o$'s formula.
\par Let the manifold $G$ be supplied with the connection.
A curve $c: [-2,2]\to G$ is called a geodesic if $\nabla \dot c(t)/dt=0$.
In view of Corollary 1.6.8 \cite{kling} there exists an open 
neighbourhood $\tilde TG$ of the submanifold $G$ of $TG$ such that 
for every $X\in \tilde TG$ the geodesic $c_X(t)$ is defined for $|t|<2$, 
where $TG$ denotes the tangent bundle. The exponential mapping
$exp^{G}: {\tilde T}G\to G$ is defined by the formula $X\mapsto c_X(1).$
The restriction $exp^G|_{\tilde TG\cap T_pG}$ will also be denoted
by $exp^{G}_p$. Then there is defined the mapping
$I(exp^{G}): I({\tilde T} G)\to I(G)$ such that for each chart
$(U,\phi )$ the mapping
$I(exp^{\phi }): Y\oplus L_{1,2}(Z,Y)\to Y\oplus L_{1,2}(Z,Y)$
is given by the following formula: 
$$I(exp^{\phi })(a^{\phi },
A^{\phi })=(a^{\phi }- \quad tr ( \Gamma ^{\phi }(A^{\phi },
A^{\phi }))/2,A^{\phi }),$$ 
where $\Gamma $ denotes the Christoffel symbol.
\par Therefore, if ${\sf R}_{x,0}(a,A)$ is a germ of diffusion processes
at a point $y=0$ of the tangent space $T_xG$, then
${\tilde exp}_x{\sf R}_{x,0}(a,A):={\sf R}_x(I(exp_x)(a,A))$ is a germ
of stochastic processes at a point $x$ of the manifold $G$.
The germs ${\tilde exp}_x{\sf R}_{x,0}(a,A)$ are called stochastic 
differentials and the It$\hat o$ bundle is called the bundle of 
stochastic differentials such that ${\sf R}_{x,0}(a,A)=:
a_xdt+A_xdw$. A section $\sf U$ of the vector bundle
$\Pi =\tau _{Y}\oplus L_{1,2}(\theta ,\tau _{Y})$ is called the 
It$\hat o$ field on the manifold $G$ and it defines a field of
stochastic differentials
${\sf R}_x(I(exp_x)(a,A))={\tilde exp}_x(a_xdt+A_xdw)$.
A random process $\xi $ has a stochastic differential defined by the 
It$\hat o$ field $\sf U:$ $d\xi (s,\omega )={\tilde exp}_{\xi (s,\omega )}
{\sf R}(a_{\xi (s,\omega )},A_{\xi (s,\omega )})$ 
if the following conditions are 
satisfied: for $\nu  _{\xi (s)}$-almost every $x\in Y$
there exists a neighbourhood $V_x$ of a point $x$ and a diffusion process
$\eta _x(t,\omega )$ belonging to the germ ${\sf R}_x(I(exp_x))(a,A)$
such that $P_{s,x} \{ \xi (t,\omega )=\eta _x(t,\omega ): 
\xi (t,\omega )\in V_x,
t\ge s \} =1$ $\nu _{\xi (s)}$-almost everywhere, where
$P_{s,x}(S):=P \{ S: \xi (s,\omega )=x \} $, $S$ is a $P$-measurable
subset of $\Omega $, $\nu _{\xi (s)}(F):= P \{ \omega :
\xi (s,\omega )\in  F \} $ (see Chapter 4 in \cite{beldal}).                     
\par If ${\sf U}(t)=(a(t), A(t))$ is a time dependent It$\hat o$ field,
then a random process $\xi (t,\omega )$ 
having for each $t\in [0,T]$ a stochastic 
differential $d\xi =exp_{\xi (t,\omega )}(a_{\xi (t,\omega )}dt+
A_{\xi (t,\omega )}dw)$
is called a stochastic differential equation on the manifold
$G$, the process $\xi (t,\omega )$ is called its solution (see Chapter 
VII in \cite{dalfo}).
As usually a flow of $\sigma $-algebras consistent with
the Wiener process $w(t,\omega )$ is a monotone set
of $\sigma $-algebras ${\sf F}_t$ such that $w(s,\omega )$ is
${\sf F}_t$-measurable for each $0\le s\le t$ 
and $w(\tau ,\omega )-w(s,\omega )$ is independent from ${\sf F}_t$
for each $\tau >s\ge t$, where ${\sf F}_s\supset {\sf F}_t$ 
for each $0\le t\le s$.
If $G$ is the manifold with the uniform atlas (see \S 2.1), 
the It$\hat o$ field
$(a,A)$ and Christoffel symbols are bounded, then
there exists the unique up to stochastic equivalence 
random evolution family $S(t,\tau )$ consistent
with the flow of $\sigma $-algebras ${\sf F}_t$
generated by the 
solution $\xi (t,\omega )$ of the stochastic differential equation 
$d\xi =exp_{\xi (t,\omega )}(a_{\xi (t,\omega )}dt+A_{\xi (t,\omega )}dw)$
on $G$, that is, $\xi (\tau ,\omega )=x$, $\xi (t,\omega )=S(t,\tau ,
\omega )x$
for each $t_0\le \tau <t<\infty $ (see Theorem 4.2.1 \cite{beldal}).

\par Address: Theoretical Department, Institute of General Physics,
\par Russian Academy of Sciences,
\par Str. Vavilov 38, Moscow, 119991 GSP-1, Russia
\end{document}